\documentclass[12pt]{amsart}
\usepackage[dvips]{epsfig}
\usepackage{graphics}
\usepackage{latexsym}
\usepackage{verbatim}
\usepackage{amsmath}
\usepackage{amsthm}
\usepackage{amssymb}
\newtheorem{theorem}{Theorem}[section]

\def\Remark{\medskip\noindent{\bf Remark: }}
\def\Remarks{\medskip\noindent{\bf Remarks: }}

\def\bul{{$\bullet$\hspace*{2mm}}}

\newcommand{\ens}[1]{\mathbb{#1}}

\newcommand{\N}{\mathbb{N}}

\newcommand{\R}{\mathbb{R}}
\newcommand{\C}{\mathbb{C}}

\def\cal{\mathcal}

\def\F{{\cal F}}

\def\var{\varepsilon}

\newcommand{\ud}{{\rm d}}
\newcommand{\nab}{\nabla_{x}}

\newcommand{\tor}{\mathbb{T}^{N}}

\let\cd\cdot

\newcommand{\dt}{\frac{\ud}{\ud t}}

\newcommand{\nabv}{\nabla_{v}}

\def\var{\varepsilon}
\def\signcm{\bigskip\bigskip\hspace{80mm}
\vbox{{\sc C. Mouhot\par\vspace{3mm}
Unit\'e de Math\'ematiques Pures et Appliqu\'ees, \'ENS Lyon\par
46 all\'ee d'Italie\par
69364 Lyon Cedex 07\par
FRANCE\par\vspace{3mm}
e-mail:} cmouhot@umpa.ens-lyon.fr }}
\def\signln{\bigskip\bigskip\hspace{80mm}
\vbox{{\sc L. Neumann\par\vspace{3mm}
Institut für Mathematik, 
Universit\"at Wien\par
Nordbergstrasse 15\par
1090 Wien\par
Austria \par\vspace{3mm}
e-mail:} Lukas.Neumann@univie.ac.at }}

\begin{document}

\title[Quantitative perturbative study of collisional kinetic models]
{Quantitative perturbative study of convergence to equilibrium for 
collisional kinetic models in the torus}

\author{Cl\'ement Mouhot, Lukas Neumann}

\hyphenation{bounda-ry rea-so-na-ble be-ha-vior pro-per-ties
cha-rac-te-ris-tic}

\begin{abstract} 
For a general class of linear collisional kinetic models in the torus, 
including in particular the linearized Boltzmann equation 
for hard spheres, the linearized Landau equation with hard and moderately soft potentials 
and the semi-classical linearized fermionic and bosonic relaxation models, 
we prove explicit coercivity estimates on the associated integro-differential operator 
for some modified Sobolev norms. 
We deduce existence of classical solutions near equilibrium 
for the full non-linear models associated, with explicit regularity 
bounds, and we obtain explicit estimates on the rate of exponential 
convergence towards equilibrium in this perturbative setting. 
The proof are based on a linear energy method which combines 
the coercivity property of the collision operator in the velocity space 
with transport effects, in order 
to deduce coercivity estimates in the whole phase space.  
\end{abstract}

\maketitle




\tableofcontents

\section{Introduction}
\setcounter{equation}{0}

In this paper we study a general class of linear inhomogeneous 
kinetic equations in the torus (including linearized 
Boltzmann, Landau, classical relaxation, semi-classical relaxation and 
Fokker-Planck equations). We then use the properties obtained on 
their evolution semi-groups to gain insight into the behaviour of the full non-linear models in the 
cases of Boltzmann, Landau and semi-classical relaxation equations.
The main tool is an estimate of the following type. 
Assume that the linear (self-adjoint) collision operator $L$ acting only on 
the velocity space is coercive for a certain norm, 
and that $L$ has a structure ``mixing part + coercive part'' 
(which shall be given a precise meaning in Subsection~\ref{subsec:model}). 
Then the integro-differential operator  $T=L-v\cd\nabla_x$ (taking into account transport effects) acting 
on the phase space of positions $x$ and velocities $v$, 
satisfies a coercivity property in $x$ and $v$, which 
implies in particular the existence of a spectral gap estimate when 
$L$ has a spectral gap in velocity (note that in general $T$ is not sectorial). 
Moreover our proof shows how to compute the constant of 
coercivity of $T$ according to the one of $L$ in the velocity space 
(and thus the rate of exponential convergence to equilibrium 
when $L$ has a spectral gap in velocity). 
Before we explain our method and results in more detail, 
let us introduce the models and problems in a precise way. 

\subsection{The Problem and its motivation}\label{Subs:PM}

In order to study the convergence to equilibrium, 
a new quantitative method in the large, 
the ``entropy-entropy-production'' method (EEP-method), has 
developed from the beginning of the 1990 decade, see 
\cite{CC92,CC94,ToVi99} in the 
spatially homogeneous setting, and \cite{DV01} in the spatially 
inhomogeneous setting. It has provided new powerful and robust tools 
in the study of relaxation towards equilibrium (see for instance  
\cite{CCG03,DV01,DV04,FNS04,NS04}), and seems to be 
the best suitable approach to deal with non-linear models in the large. 
However, it can be seen from these references that the method, while 
very robust with respect to nonlinearities and able to deal with external potentials at 
no increased difficulties, has two major shortcomings. First it relies on uniform in time 
estimates of regularity of the solution, that are usually hard to establish. Second, 
it seems to fail to give the optimal rate of convergence, in particular when $L$ has 
a spectral gap in velocity it seems to fail to give exponential rate of 
convergence to equilibrium. 

In this article we provide a linear energy method in order to overcome these problems for linear models, 
or for non-linear models in the perturbative setting (this method 
is explicit in terms of the explicit coercivity estimates on the linearized 
collision operators in the velocity space, see~\cite{BaMo,Mcoerc,MS**}).  
The difference when compared to the EEP-method is that the convergence as well as 
the uniform in time regularity bounds are obtained in a single step. However this approach 
is linear and therefore limited to the perturbative setting near equilibrium for the non-linear 
models. A similar approach has previously been used in Guo's papers 
\cite{Guo:BE:03,Guo:VMB:03,Guo04,SG}. 
Also another approach has been developed in~\cite{HeNi} for the Fokker-Planck 
model in a confining external field (see also~\cite{He:VPFP}), and H\'erau recently generalized the method to the linear 
relaxation model in~\cite{He}. Similar results for the Fokker-Planck equation 
were recovered in~\cite{Vi:hypo} by an energy method. 
The results we present here are partly included in these references, 
however the proof is simpler and more explicit. Our viewpoint unifies 
previous scattered results and explores new situations such as semi-classical 
relaxation models. The works quoted so far are the starting point of this paper. 
In particular one of the key ideas of our proof, namely looking for the time derivative 
of a ``mixed term'' of the form $\int \nabla_x h \cdot \nabla_v h$, is inspired 
from~\cite{Vi:hypo} where it first appeared in this explicit form (see also~\cite{Vi:hypofd}).

Most of these articles deal with the Boltzmann equation that we address in section~\ref{linbol}.
Our main abstract theorem applies to the Boltzmann linearized collision 
operator for hard spheres (and hard potentials with Grad's cutoff assumption): 
smooth solutions with explicit regularity bounds and rate of convergence to 
equilibrium are constructed near equilibrium. 
A similar study has been performed 
in a non-constructive way for the Vlasov-Poisson-Boltzmann system near equilibrium in \cite{Guo:PB:02} 
and for the Landau equation in \cite{Guo:LE:02}. In \cite{Guo:BE:03} the same method 
was also applied to the Boltzmann equation for cutoff soft potentials. 

Guo's argument relies, roughly speaking, on the coercivity of 
the linearized Boltzmann operator ``in the mean'' -- {\em i.e.}, if 
integrated over a time interval -- for small perturbations of the Maxwellian. 
Due to this averaging in time precise rates of convergence are quite hard to 
deduce from the results in \cite{Guo:PB:02,Guo:LE:02}. 
Later Guo refined his approach in \cite{Guo:VMB:03} to meet 
the needs of application to the Vlasov-Maxwell-Boltzmann system. 
The new idea is to use a norm that includes temporal derivatives. 
In this norm the fields are almost controlled by the deviation of the 
distribution from the Maxwellian -- to be more precise they ``loose'' one derivative. 
This leads to \emph{instantaneous} coercivity and thus global perturbative solutions 
for the Vlasov-Maxwell-Boltzmann system as well as exponential decay to 
equilibrium in the case of Vlasov-Poisson-Boltzmann, in these norms. 
The proof, while being very complicated, is constructive. 
More recently a farther refined method based also on norms including temporal 
derivatives has been used to show almost exponential decay to equilibrium for various of 
the mentioned models and the relativistic Landau Maxwell system in \cite{SG}.

The EEP-method has been applied to the Boltzmann equation in \cite{DV04} -- 
again assuming uniform in time regularity bounds on the solution -- 
but for a large class of collision operators and most important -- without perturbative assumptions. 
All these works have been carried out on bounded domains, essentially the torus apart from \cite{DV04} 
where various types of boundary conditions are considered. For completeness we also refer to a 
recent preprint dealing with the problem in whole space \cite{Guo04} and the article \cite{LYY04} 
for a related non-linear energy method. 

Despite this vast amount of recent literature on the 
problem our proof has two advantages. First it is simpler than the 
ones in the articles quoted above, 
mostly due the fact that we take advantage of the mixing 
properties of the collision operator. 
Second it separates in a very clear way between linear effects in these energy 
estimates (transport + linear collision), which are expressed in a coercivity estimate 
on this linear part, and the problems arising from the small remaining 
bilinear part when considering solutions near the equilibrium. 
We are able to derive exponential convergence to equilibrium 
without resorting to norms including time derivatives and in a purely instantaneous manner.  

\subsection{The models}\label{subsec:model}

We will study initial value problems for equations of the form 
\begin{equation}\label{eqgenNL}
\partial_t f + v \cd \nab f = Q(f), \qquad t \ge 0, \ x \in \Omega, \ v \in \R^N,
\end{equation} 
where $N \ge 1$ denotes the dimension of space. 
In this equation, $f$ denotes the distribution 
of particles in phase space, therefore it is a time-dependent non-negative 
$L^1(\Omega\times\mathbb{R}^N)$ function. The operator $Q$ models the collisional 
interactions of particles (either binary or between 
particles and a surrounding medium). It is local in $x$ and $t$ 
and it depends on the particular model of interaction chosen. 
The term $v\cd \nab f$ corresponds to the free flow of particles. 
For the spatial domain $\Omega$, we shall consider here 
the periodic case, that is $\Omega=\tor$. Hence $f=f(t,x,v)$ satisfies 
initial conditions $f(t=0,x,v)=f_0(x,v)$ and periodic boundary conditions. 

One defines the global mass, momentum and energy of the solution $f$ as 
  \[ \int_{\Omega \times \R^N} f \, \ud x \, \ud v, \qquad \int_{\Omega \times \R^N} f \, v \, \ud x \, \ud v, \qquad 
     \int_{\Omega \times \R^N} f \, |v|^2 \, \ud x \, \ud v. \]
Depending on the model chosen for the collision operator $Q$, one has 
preservation along time of part or all of these quantities. 

Under very general assumptions, equations of class~\eqref{eqgenNL} 
admit a unique equilibrium in the torus, 
which we shall denote by $f_{\infty}$, and 
which is independent of $t,x$ (this is trivial for relaxation models admitting 
only mass conservation, and for models admitting conservation of 
mass, momentum and energy, such as Boltzmann and Landau equations, 
this is shown easily inspiring from the arguments in~\cite{DV04} and~\cite{De90}). 

Then one can consider the linearization 
around this equilibrium. In order to reduce to a Hilbert space setting, 
one usually considers perturbations of the form 
$f=f_\infty + f_\infty^{1/2} h$. Discarding the bilinear 
term, it yields the following linearized equation on $h$   
  \begin{equation}\label{eqgenL}
  \partial_t h +v\cd\nab h=L(h),
  \end{equation}
where $L$ depends on the precise form of the collision operator $Q$. 
The unknown $h$ belongs to $L^2 (\tor \times \R^N)$. 


Now let us give a general framework for linear collisional kinetic models.  
We shall denote as usual by $L^2$ the Lebesgue space of square integrable functions, 
and, for $k \in \N$, $H^k$ the Sobolev space of $L^2$ functions with square integrable 
derivatives up to order $k$. When needed we shall indicate as a subscript the 
variables ($x$ or $v$) these functional spaces refer to. For the variable 
$x$, Sobolev spaces refer implicitly to {\em periodic} Sobolev spaces. 
When no subscript is used, these functional spaces always refer to 
$x$ and $v$.   
\medskip

\begin{itemize}
\item[{\bf H1.}] (\underline{Structure}) We consider 
a linear collision operator $L$ on $L^2 = L^2 (\tor \times \R^N)$ 
which a closed and self-adjoint operator on $L^2 _{v}$  
and local in $t,x$. 

We assume that it writes 
  \begin{equation}\label{eq:Lgen}
  L = K - \Lambda 
  \end{equation}
where {\bf $\Lambda$ is a coercive operator} in the sense: 
there is a norm $\|\cdot\|_{\Lambda_v}$ on $\mathbb{R}^N$ (the space of velocities), 
such that 
  \[ \nu^\Lambda _0 \, \|h\|^2 _{L^2 _v} \le 
      \nu^\Lambda _1 \|h\|^2 _{\Lambda_v} 
       \le \langle \Lambda(h), h \rangle_{L^2 _v} 
       \le \nu^\Lambda _2 \,  \|h\|^2 _{\Lambda_v} \]
and also 
  \[ \langle \nabla_v \Lambda(h), \nabla_v h \rangle_{L^2 _v} \ge 
      \nu^\Lambda _3 \, \| \nabla_v h\|^2 _{\Lambda_v} - \nu^\Lambda _4 \, \|h\|^2 _{L^2 _v} \]
for some constants $\nu^\Lambda _0, \nu^\Lambda _1, \nu^\Lambda _2, \nu^\Lambda _3, \nu^\Lambda _4 >0$. 

Moreover we assume that $L$ satisfies (for some constant $C^L >0$) 
  \[ \langle L(h), g \rangle_{L^2 _v} \le C^L \,  \|h\| _{\Lambda_v} \, \|g\|_{\Lambda_v}. \]
To shorten the notation we also introduce the norm 
\begin{equation*}
\|\cdot\|_{\Lambda}:=\left\|\|\cdot\|_{\Lambda_v}\right\|_{L^2_x}\ .
\end{equation*}
\medskip

\item[{\bf H2.}] (\underline{Mixing property in velocity}) 
We assume that {\bf $K$ has a regularizing effect} in the following sense: for 
any $\delta >0$, there is some explicit $C(\delta)$ such that 
for any $h \in H^1 _v$, 
  \begin{equation}\label{eq:Kgen}
  \langle \nabla_v K(h), \nabla_v h \rangle_{L^2 _v} 
  \le C(\delta) \, \|h\|^2 _{L^2 _{v}} + \delta \, \|\nabla_{v}h\|^2 _{L^2 _{v}} \ .
  \end{equation}
\medskip

\item[{\bf H3.}] (\underline{Relaxation to the local equilibrium}) 
We assume that $L$, as an operator on $L^2 _v$, has a finite dimensional kernel  
  \[ N(L) = \mbox{Span} \left\{ \varphi_1, \dots, \varphi_n  \right\} \]  
and we denote by $\Pi_l$ the orthogonal projection on $N(L)$ in $L^2 _v$.  
We make the following {\bf local coercivity assumption}: 
there is $\lambda>0$ such that  
  \[ \langle L(h), h\rangle_{L^2 _v} \le - \lambda \, \|h - \Pi_l(h)\|_{\Lambda_v} ^2. \]
This together with {\bf H1} implies in particular that $L$ is non-positive  
and has a spectral gap in $L^2_v$, whose size is bounded from below by 
  \begin{equation*}
  \lambda_L=\left( \frac{\nu_0^{\Lambda}}{\nu_1^{\Lambda}} \right) \lambda>0.
  \end{equation*}
\end{itemize}
\medskip

In the structure assumption {\bf H1}, $K$ shall typically stand  
for a multiplicative operator by a function $\nu$ 
(usually called the {\em collision frequency}) in the case 
of short-range interactions, or some diffusion operator in the 
case of long-range interactions. The norm $\|\cdot\|_\Lambda$ that 
we define can also loosely be seen as the norm of the graph of $\Lambda^{1/2}$.  
The coercivity of $\Lambda$ is linked with the coercivity property of the whole linearized collision 
operator $L$ in {\bf H3} (where the word ``local'' refers to the position $x$ and
the fact that it pushes the dynamic towards local 
equilibrium). This coercivity property is crucial to ensure exponential decay towards 
equilibrium, even in the homogeneous setting. In the case where it is 
weakened, which happens for interactions with ``weak collision effect'' 
such as soft potentials (see for instance in \cite{Mcoerc}), one expects convergence rates 
of the form $e^{-t^\tau}$ with $\tau<1$ (see the results in~\cite{SG:soft}). 

The kernel of $L$ in $L^2 _v$, which is composed of functions 
belonging to $N(L)$ for any $x$, corresponds to the manifold of local equilibria for the linearized kinetic models. Therefore when the $x$ variable is added, 
$\Pi_l$ is the projection on the ``fluid part'', and $(\mbox{Id}-\Pi_l)$ is the projection 
on the ``kinetic'' part. It is defined by 
  \[ \Pi_l(h) = \sum_{i=1} ^n \left( \int_{\R^N} h \, \varphi_i \, \ud v \right) \, \varphi_i. \]

It is trivial that any 
local equilibrium uniform in space is indeed a global equilibrium. Since $L$ is self-adjoint, 
the $\varphi_1, \dots, \varphi_n$ belong to the kernel of its adjoint 
$L^*$, and thus by integrating in $x$ and $v$ we get 
  \[ \forall \, i = 1, \dots , n, \qquad \frac{\ud}{\ud t} 
         \int_{\tor \times \R^N} h \, \varphi_i \, \ud x \, \ud v = 0. \]
Hence we denote  
  \[ \Pi_g(h) = \sum_{i=1} ^n \left( \int_{\tor \times \R^N} h \, \varphi_i \, \ud x \, \ud v \right) \, \varphi_i\ \]
which is time and space independent (this can be shown easily to be the 
orthogonal projection on $N(T)$ in $L^2_{x,v}$, directly or using the 
coercivity property of Theorem~\ref{theo:T}).   
\smallskip

A detailed study of physical models satisfying assumptions {\bf H1-H2-H3} shall be 
given in Section~\ref{sec:phys}. 

\subsection{Main results}

First we state the main result on the coervivity estimates on $T$ and the 
consequence on its evolution semi-group. 

  \begin{theorem}\label{theo:T}
  Let $L$ be a linear operator on $L^2 _{x,v}$ satisfying assumptions 
  {\bf H1-H2-H3}. Then $T=L-v \cdot \nabla_x$ generates a strongly continuous evolution 
  semi-group $e^{tT}$ on $H^1 _{x,v}$, which satisfies 
    \[ \| e^{tT} (\mbox{{\em Id}} - \Pi_g) \|_{H^1 _{x,v}} \le C_T \, \exp \left[ -\tau_T t \right] \]
  for some explicit constants $C_T$, $\tau_T >0$ depending 
  only on the constants appearing in assumptions {\bf H1-H2-H3}. 
  More precisely
    \[ \forall \, h \in H^1, \quad 
       \langle T h , h \rangle_{\cal{H}^1} \le - C_T ' \, \left( \| h - \Pi_g (h) \|^2 _{\Lambda}
                                       + \| \nabla_{x,v} (h - \Pi_g (h)) \|^2 _{\Lambda} \right) \]
  for some (explicit) Hilbert norm $\cal{H}^1$ equivalent to the $H^1$ norm, 
  and some explicit constant $C_T '>0$. 
  \end{theorem}

\Remarks 

1. Note that under slightly strengthened assumptions, a similar result is proved in 
$H^k$ in Theorem~\ref{theo:T:k}. 
\smallskip

2. The method does not rely on an abstract result from 
spectral theory such as Weyl's theorem, like Ukai's proof of the existence of 
a spectral gap for the Boltzmann equation for hard spheres in~\cite{Uka74}. 
Hence we do not need the compactness property of $K$, although we require 
a regularizing property on $K$ which is strongly related 
(see the discussion in~\cite{Meepts}).  
Our method can be seen as a quantitative version of Ukai's result 
(in the case of the linearized Boltzmann equation). In particular 
it shows that apart from $0$, the spectrum of $T$ is included in 
  \[ \{\xi \in \C \, ; \, \mbox{Re}(\xi) \le - \tau_T \}. \] 
The abstract setting emphasizes what is effectively required 
from the linearized collision operator to 
deduce exponential convergence, and it allows to apply the method 
to other models as well. Since for the linearized Boltzmann equation 
it was proved in~\cite{Uka74} that $T$ is not sectorial (its essential 
spectrum is given by a half-plane), our work can also be seen as a method 
to prove exponential decay of the semi-group for a class of 
non-sectorial operators, which is in general quite tricky. In the case of 
Fokker-Planck type operators, other methods have been developed 
in~\cite{HeNi,Vi:hypo} to solve this question. 
\smallskip

3. In order to obtain a completely quantitative result of 
convergence to equilibrium, one has to get estimates on 
the constant $\lambda$ in assumption {\bf H3}. This question had 
remained open for a long time for important physical models such as the 
linearized Boltzmann collision for hard spheres or the Landau collision 
operators for hard and moderately soft potentials. It has been solved 
recently in the works~\cite{BaMo,Mcoerc,MS**}. Therefore this theorem 
allows to compute rates of convergence explicitly for all the models 
we consider in Section~\ref{sec:phys} (except semi-classical relaxation for 
bosons). 
\medskip

In the two following theorems we define $k_0$ the smallest integer 
such that $E(k_0/2) >N/2$ (where $E$ denotes 
the integer value). Our second main result is for the nonlinear 
Boltzmann and Landau models. 

\begin{theorem}\label{theo:NL}
Let us consider either the Boltzmann equation~\eqref{eq:Bol} for hard 
spheres or hard potentials with cutoff or the Landau equation~\eqref{eq:Lan} 
with $\gamma \ge -2$, in the torus.  
Let $0 \le f_0$ be an initial datum with finite mass and energy, and  
we denote by $f_\infty$ the unique equilibrium associated to $f_0$. 
We suppose that the initial datum satisfys
  \[ \| f_{\infty}^{-1/2} \, ( f_0 - f_\infty) \|_{H^k} \le \varepsilon \]
for some $k \ge k_0$ and some $0< \varepsilon \le \varepsilon_0$ 
where $\varepsilon_0$ depends explicitly on the collision operator.  

Then there exists a unique global 
solution $0 \le f=f(t,x,v)\in\mathcal{C}([0,\infty[,H^k)$ of 
the initial value problem \eqref{eqgenNL}, such that 
  \[ \forall \, t \ge 0, \qquad \| f_\infty ^{-1/2} ( f(t,\cd,\cd) - f_\infty) \|_{H^k} \le 
      C \, \exp \left[ -\tau t \right] \]
for some explicit constants $C, \tau >0$. 

The conclusion still holds true when a repulsive self-consistent 
Poisson potential is added (without smallness condition on the intensity of 
the self-consistent interaction), in the case of the Boltzmann equation for hard spheres, or 
the Landau equation with $\gamma \ge -1$.  
\end{theorem}
The next theorem deals with the semi-classical relaxation models.
\begin{theorem}\label{theo:NLSC}
Consider the semi-classical relaxation equation~\eqref{eqfsc} in the torus for fermions ($\epsilon=1$) or 
bosons ($\epsilon=-1$), with an initial datum $0 \le f_0 \in L^1$. 
Let the equilibrium distribution be given by 
\begin{equation*}
f_\infty=\frac{\kappa_\infty \, \mathcal{M}}{1+\epsilon\, \kappa_\infty \, \mathcal{M}},
\end{equation*}
with $1+\epsilon \, \kappa_\infty \, \mathcal{M}>0$, where $\mathcal{M}$ is 
the normalized Maxwellian 
 \[ {\cal M}(v) = \frac{e^{-|v|^2/2}}{(2\pi)^{N/2}} \]
and $\kappa_\infty$ is defined by mass conservation.
Let $k \ge k_0$ and let the initial datum satisfy
  \[ \left\| \,(1+\epsilon \, \kappa_\infty \, \mathcal{M}) 
     \left(\kappa_\infty \, \mathcal{M}\right)^{-1/2} ( f_0 - f_\infty ) \right\|_{H^k} \le \varepsilon \]
for some $0< \varepsilon \le \varepsilon_0$ 
where $\varepsilon_0$ depends explicitly on the collision operator.  

Then there exists a unique global 
solution $0 \le f=f(t,x,v)\in\mathcal{C}([0,\infty[,H^k)$ of 
the initial value problem \eqref{eqgenNL}, such that 
  \[ \forall \, t \ge 0, \qquad \left\| (1+\epsilon \, \kappa_\infty \, \mathcal{M})
     \left(\kappa_\infty \, \mathcal{M}\right)^{-1/2} ( f(t,\cd,\cd) - f_\infty) \right\|_{H^k} \le 
      C \, \exp \left[ -\tau t \right] \]
for some explicit constants $C, \tau>0$.  
\end{theorem}
\Remark
The condition on the form of the equilibrium distribution is 
in fact trivially fulfilled for the fermionic case.
In the bosonic case it is a condition of smallness on $\kappa_\infty$ 
and thus on the initial mass. Indeed it is equivalent to impose that 
the mass of the initial datum is small enough such that no condensation 
occurs. This is not for technical reasons but necessary to ensure 
exponential convergence as can be seen from the detailed study 
of the asymptotics in the spatially homogeneous case in \cite{EMV}.
\medskip

\subsection{Outline of the article}
The article is structured as follows.
In Section~\ref{sec:main} we give the proof of Theorem~\ref{theo:T}, which turns 
out to be very short and simple. Then in Section~\ref{sec:ext} we 
expose several extensions of the method. In particular we show how to generalize Theorem~\ref{theo:T} 
to higher-order Sobolev norms, and how to include a weak external field or a self-consistent 
Poisson potential in the study. 
Section~\ref{sec:NL} is devoted to the application of the previous study to 
the genuine non-linear problems of the form~\eqref{eqgenNL} 
near equilibrium, and we prove the abstract Theorem~\ref{theo:pert}. 
Finally in Section~\ref{sec:phys}, we discuss the general 
assumptions of Theorem~\ref{theo:T} and Theorem~\ref{theo:pert} for 
an extensive list of physical models: 
classical or semi-classical relaxation, Boltzmann equation 
for hard spheres or hard potentials with cutoff, 
Landau equation for hard potentials or moderately soft potentials, linear 
Fokker-Planck equation. Then Theorem~\ref{theo:NL} follows from 
this study together with the Theorem~\ref{theo:pert} (we also comment on the marginal  
differences between the proofs of Theorem~\ref{theo:NL} and Theorem~\ref{theo:NLSC}).

\section{Proof of Theorem~\ref{theo:T}}\label{sec:main}
\setcounter{equation}{0}

We divide the proof into several steps. 
\medskip

\noindent
{\bf Step 1.} Since by assumption {\bf H3}, the operator $L$ is non-positive, the $L^2 _{x,v}$ norm is 
decreasing along the flow and it is straightforward to deduce that $T$ generates a 
strongly continuous contraction evolution semi-group on $L^2 _{x,v}$ (see~\cite{Kato} 
for instance). 
In order to estimate the semi-group in $H^1 _{x,v}$, let us consider 
$h_0 \in H^1 _{x,v} \cap \mbox{Dom}(T) \cap N(T)^\bot$ and 
$h=h(t,x,v)$ the associated solution of the equation $\partial_t h = T(h)$. 
Now we study the evolution of the $H^1 _{x,v}$ norm of $h$. 
\medskip

\noindent
{\bf Step 2.} 
We estimate the time evolution of the $L^2$ norm of $h$. Using the skew symmetry 
of the transport part together with assumption {\bf H3}, we find immediately
  \begin{equation}\label{coercL^2}
  \dt \| h \|^2 _{L^2} \le - 2 \, \lambda \,  \| h - \Pi_l h \|^2 _\Lambda.
  \end{equation}
\medskip

\noindent
{\bf Step 3.} 
We estimate the time derivative of the gradients in $x$ and $v$. 
\smallskip

\bul For the gradient with respect to $x$ we obtain, thanks to assumption {\bf H3},  
  \begin{equation}\label{eqx}
  \dt \|\nabla_x h\|^2 _{L^2} \le - 2 \, \lambda \, \|\nabla_x h - \Pi_l(\nabla_x h) \|^2 _\Lambda.
  \end{equation}

\bul For the gradient with respect to $v$ we get 
  \[ \dt \left\|\nabla_v h\right\|^2 _{L^2} 
     = 2 \, \langle \nabla_v K(h), \nabla_v h \rangle _{L^2}
                                     -2 \, \langle \nabla_v \Lambda(h),\nabla_v h \rangle_{L^2} 
                                     - 2 \, \langle \nabla_x h, \nabla_v h \rangle_{L^2}. \]
Using assumption {\bf H1} we have 
  \[ -2 \, \langle \nabla_v \Lambda(h),\nabla_v h \rangle_{L^2} 
       \le - 2 \, \nu^\Lambda _3 \, \|\nabla_v h\|^2 _\Lambda 
           +  2 \, \nu^\Lambda _4 \, \|h\|^2 _{L^2}. \]
Assumption {\bf H2} yields 
\begin{multline*}
   2 \, \langle \nabla_v K(h), \nabla_v h \rangle_{L^2}  
          \le \frac{\nu^\Lambda _1}{\nu^\Lambda _0} \, \frac{\nu^\Lambda _3}{2} \, \left\|\nabla_v h\right\|^2 _{L^2}
           + 2 \, C \, \left(\frac{\nu^\Lambda _1}{\nu^\Lambda _0}\frac{\nu^\Lambda _3}{4}\right) \, 
           \left\| h \right\|^2_{L^2} \leq\\
           \frac{\nu_3 ^\Lambda}{2} \, \|\nabla_v h\|_{\Lambda}^2
           + 2 \, C \, \left(\frac{\nu^\Lambda _1\nu^\Lambda _3}{4\nu^\Lambda _0}\right) 
           \left\| h \right\|^2_{L^2}. 
\end{multline*}            
Furthermore we split
  \begin{multline*}
  2 \, \langle \nabla_x h, \nabla_v h \rangle _{L^2} 
         \le \frac{\nu^\Lambda _1}{\nu^\Lambda _0}\frac{\nu^\Lambda _3}{2} 
         \, \left\|\nabla_v h\right\|^2 _{L^2} 
         + \frac{\nu^\Lambda _0}{\nu^\Lambda _1}\frac{2}{\nu^\Lambda _3} \, 
        \left\|\nabla_x h\right\|^2 _{L^2} \\ \leq \frac{\nu^\Lambda _3}{2} \, \left\|\nabla_v h\right\|^2 _{\Lambda}
         + \frac{2\nu^\Lambda _0}{\nu^\Lambda _1\nu^\Lambda _3}\, \left\|\nabla_x h\right\|^2 _{L^2}.
  \end{multline*}
Using the last three inequalities we have
  \[ \dt \|\nabla_v h\|^2 _{L^2} \le 
             \left[ 2C\left(\frac{\nu^\Lambda _1\nu^\Lambda _3}{4\nu^\Lambda _0}\right) 
     +  2 \, \nu^\Lambda _4  \right] \, \|h\|_{L^2}^2
              + \frac{2\nu^\Lambda _0}{\nu^\Lambda _1\nu^\Lambda _3} \, \left\|\nabla_x h\right\|^2 _{L^2} 
               - \nu^\Lambda _3 \, \left\|\nabla_v h\right\|^2 _{\Lambda}. \]
Now we write
  \[ \left\|h\right\|^2 _{L^2} \le 2 \, \left\|h-\Pi_l(h) \right\|^2 _{L^2} 
                            + 2 \, \left\| \Pi_l (h) \right\|_{L^2} ^2.\]
Since $\Pi_g(h) =0$ we deduce that $\Pi_l(h)$ has zero mean on the torus, 
and Poincar\'e's inequality in the torus yields 
(for a constant $C_P$ only depending on the dimension $N$) 
  \[ \| \Pi_l (h) \|^2_{L^2} \le C_P \, \|\Pi_l (\nabla_x h) \|^2_{L^2} \le C_P \, \| \nabla_x h \|^2_{L^2}, \]
and thus we get for some explicit constants $C_1,C_2>0$
  \begin{equation}\label{main:v}
   \dt \|\nabla_v h\|^2 _{L^2} \le 
          C_1 \, \|h - \Pi_l(h) \|_{\Lambda} ^2 + C_2 \, \|\nabla_x h\|^2 _{L^2} 
    - \nu^\Lambda _3 \,  \|\nabla_v h\|^2 _\Lambda. 
  \end{equation}                             

\bul For the mixed term we have 
  \[ \dt \langle \nabla_x h, \nabla_v h \rangle _{L^2} = 
               - \|\nabla_x h\|^2 _{L^2} + 2 \, \langle \nabla_x L(h), \nabla_v h \rangle _{L^2}. \]
Then we write (using assumption {\bf H1} and Cauchy-Schwarz's inequality in $x$)
  \begin{multline*}
  2 \, \langle \nabla_x L(h), \nabla_v h \rangle _{L^2} = 
        2 \, \langle L \left( \nabla_x h - \Pi_l(\nabla_x h)\right), \nabla_v h \rangle _{L^2} \\ 
        \le 2 \, C^L \, \left\| \nabla_x h - \Pi_l(\nabla_x h)\right\| _\Lambda
              \, \left\| \nabla_v h \right\| _\Lambda \\
         \le C^L \, \eta \, \left\| \nabla _x h - \Pi_l(\nabla_x h)  \right\|^2 _\Lambda
              + C^L \, \eta^{-1} \, \left\| \nabla_v h \right\|^2 _\Lambda 
  \end{multline*}
for any $\eta >0$. Hence we obtain
  \begin{equation}\label{eqxv}
  \dt \langle \nabla_x h, \nabla_v h \rangle \le - \left\|\nabla_x h\right\|^2 _{L^2} 
             + C^L \, \eta \, \left\| \nabla_x h - \Pi_l(\nabla_x h) \right\|^2 _\Lambda 
              + C^L \, \eta^{-1} \, \left\| \nabla_v h \right\|^2 _\Lambda.
  \end{equation} 
\medskip

\noindent
{\bf Step 4.} 
Now it remains to combine equations~\eqref{coercL^2}, 
\eqref{eqx}, \eqref{main:v} and~\eqref{eqxv}: we pick $A, \alpha, \beta, \gamma >0$ and compute 
  \begin{multline*}
  \dt \Big[ A \, \|h\|^2 _{L^2} + \alpha \, \|\nabla_x h\|^2 _{L^2}
                         + \beta \, \|\nabla_v h\|^2 _{L^2}
                         + \gamma \, \langle \nabla_x h, \nabla_v h \rangle _{L^2} \Big] \\ 
  \le (\beta C_1 - 2 A \lambda) \, \|h - \Pi_l(h) \|_{\Lambda} ^2  
     + \left(\eta \gamma C^L - 2 \alpha \lambda \right) \, \| \nabla_x h - \Pi_l(\nabla_x h) \|^2 _\Lambda \\
      + \left( \eta^{-1} \gamma C^L - \beta \nu^\Lambda _3 \right) \,  \|\nabla_v h\|^2 _\Lambda
      + \left( C_2 \beta - \gamma \right) \, \|\nabla_x h\|^2 _{L^2}. 
  \end{multline*}

For a given $\beta$ with $\nu_3^\Lambda \beta>1$, we first fix $A$ big enough such that 
  \[ (\beta C_1 - 2 A \lambda) \le -1, \]
then $\gamma$ big enough such that 
  \[ \left( C_2 \beta - \gamma \right)  \le -1, \]
then $\eta$ big enough such that 
  \[ \left( \eta^{-1} \gamma C^L - \beta \nu^\Lambda _3 \right)  \le  - 1, \]
then $\alpha$ big enough such that 
  \[  \left(\eta \gamma C^L - 2 \alpha \lambda \right) \le -1 \]
and such that $\gamma^2 < \alpha \beta$ and $\alpha \ge \beta$. 
For this choice we obtain 
  \begin{multline*}
  \dt \Big[ A \, \|h\|^2 _{L^2} + \alpha \, \|\nabla_x h\|^2 _{L^2}
                         + \beta \, \|\nabla_v h\|^2 _{L^2}
                         + \gamma \, \langle \nabla_x h, \nabla_v h \rangle \Big] \\
      \le 
         - \Big[ \|\nabla_x h\|^2 _{L^2} + \|\nabla_v h\|^2 _{\Lambda} 
                  + \|h - \Pi_l(h)\|^2 _{\Lambda} + \|\nabla_x h - \Pi_l(\nabla_x h) \|_{\Lambda} \Big] .
  \end{multline*}
The function 
  \[ {\cal F}(t) = \Big[ A \, \|h\|^2 _{L^2} +  \alpha \, \|\nabla_x h\|^2 _{L^2} 
                         + \beta \, \|\nabla_v h\|^2 _{L^2} 
                         + \gamma \, \langle \nabla_x h, \nabla_v h \rangle \Big] \]
satisfies (remember that $\alpha\geq\beta\,$)
  \begin{multline*}
  A \, \|h\|^2 _{L^2} + 
    (\beta/2)\, \left[ \|\nabla_x h\|^2 _{L^2} + \|\nabla_v h\|^2 _{L^2} \right] \le {\cal F}(t) \\ \le 
       A \, \|h\|^2 _{L^2} +  (3\alpha/2) \, \left[ \|\nabla_x h\|^2 _{L^2} + \|\nabla_v h\|^2 _{L^2} \right]. 
  \end{multline*}
Moreover, since $\Pi_l(h)$ has zero mean on the torus, we have by Poincar\'e's inequality:
  \[ \|h\|^2 _{\Lambda} \le 2 \, \|h - \Pi_l(h)\|^2 _{\Lambda} + 2 \, \|\Pi_l(h)\|^2 _{\Lambda} 
             \le C \, \left( \|h - \Pi_l(h)\|^2 _{\Lambda} + 1/2\|\nabla_x h\|^2 _{L^2} \right) \]
for some explicit constant $C>0$, and similarly
  \[ \|\nabla_x h\|^2 _\Lambda 
             \le C' \, \left( \|\nabla_x h - \Pi_l(\nabla_x h)\|^2 _{\Lambda} + 1/2\|\nabla_x h\|^2 _{L^2} \right) \]
for some explicit constant $C'>0$. Hence we deduce that 
  \[ \dt \cal{F}(t) \le - K \, \left( \|h\|^2 _\Lambda + \|\nabla_{x,v} h \|^2_\Lambda \right)  \]
for some explicit $K>0$, and that $\cal{F}(t)$ is equivalent to the square of the $H^1$ norm of $h$. 
We define the norm $\mathcal{H}^1$ by
\begin{equation*}
\|\cd\|_{\mathcal{H}^1} =\Big\{A\|\cd\|_{L^2}^2+\alpha\|\nabla_x\cd\|_{L^2}^2 
+\beta\|\nabla_v\cd\|_{L^2}^2+\gamma\left<\nabla_x \cd,\nabla_v\cd\right>_{L^2}\Big\}^{1/2}.
\end{equation*}
This concludes the proof. 
%

\section{Some generalizations}\label{sec:ext}
\setcounter{equation}{0}

\subsection{Higher-order Sobolev spaces}    

In this section we show how to extend the previous method 
to higher-order Sobolev spaces. 
Let us first introduce some notations about multi-indices. 
For a multi-index $j$ in $\N^N$, we shall denote by $c_i(j)$ the value of the 
$i$-th coordinate of $j$ ($i=1,\dots, N$), and by $|j|$ the $l^1$ norm 
of the multi-index, that is $|j|=\sum_{i=1} ^N c_i(j)$. We also denote by 
$\delta_{i_0}$ the multi-index such that $c_i(\delta_{i_0})=0$ for 
$i \not = i_{0}$ and $c_{i_0}(\delta_{i_0})=1$. 
Finally for two multi-indices $j$ and $l$ in 
$\N^N$ we set the 
shorthand $\partial^j_l=\partial/\partial v_j \, \partial /\partial x_l$.

In order to treat the higher-order derivatives we shall strengthen 
assumptions {\bf H1} and {\bf H2} into 
\medskip

\begin{itemize}
\item[{\bf H1'.}] 
We assume {\bf H1}. Moreover we assume that for any $k \ge 1$, 
for any multi-indexes $j$ and $l$ such that $k=|j|+|l|$ and $|j| \ge 1$, we have 
  \[ \langle \partial^j_l \Lambda(h), \partial^j_l h \rangle_{L^2} \ge 
      \nu^\Lambda _5 \, \|\partial^j_l h\|^2 _\Lambda - \nu^\Lambda _6 \, \|h\|^2 _{H^{k-1}} \]
for some constants $\nu^\Lambda _5, \nu^\Lambda _6 >0$. 
\medskip
\item[{\bf H2'.}] 
We assume that {\bf $K$ has a regularizing effect} in the sense: 
for any $k \ge 1$, for any multi-indexes $j$ and $l$ such that $k=|j|+|l|$ and $|j| \ge 1$, for 
any $\delta >0$, there is some explicit $C(\delta)$ such that 
  \begin{equation}\label{eq:Kgen:plus}
  \langle \partial^j_l K(h), \partial^j_l h \rangle_{L^2} 
  \le C(\delta) \, \| h\|^2 _{H^{k-1}} + \delta \, \| \partial^j_l h\|^2 _{L^2} \ .
  \end{equation}
\end{itemize}
\medskip

Again these strengthened assumptions are satisfied by the physical 
models we discussed in the introduction, as we check in Section~\ref{sec:phys}.
Now we can formulate the coercivity estimate on $T$ and the consequence on its semi-group.
 \begin{theorem}\label{theo:T:k}
 Let $L$ be an operator on $L^2$ satisfying assumptions 
 {\bf H1'-H2'-H3} and $k$ be an integer number. 
 Then $T=L-v \cdot \nabla_x$ generates a strongly continuous evolution 
 semi-group $e^{tT}$ on $H^k$, which satisfies 
   \[ \| e^{tT} (\mbox{{\em Id}} - \Pi_g) \|_{H^k} \le C_T \, \exp \left[ -\tau_T t \right] \]
 for some explicit constants $C_T$, $\tau_T >0$ depending 
 only on the constants appearing in {\bf H1'-H2'-H3}. 
 More precisely
    \[ \forall \, h \in H^k, \quad 
       \langle T h , h \rangle_{\cal{H}^k} \le - C_T ' \, \left( 
              \sum_{|j|+|l| \le k} \| \partial^j _l (h - \Pi_g (h)) \|^2 _{\Lambda} \right)  \]
  for some (explicit) Hilbert norm $\cal{H}^k$ equivalent to the $H^k$ norm, 
  and some explicit constant $C_T '>0$. 
 \end{theorem}

\begin{proof}[Proof of Theorem~\ref{theo:T:k}] 
We write $\dot H^k$ for the homogeneous Sobolev semi-norm, {\it i.e.}, 
  \[ \|h\|_{\dot H^k} ^2 =\sum_{|j|+|l|=k}\|\partial^j _l h\|_{L^2} ^2. \]
Again we pick $h_0\in H^k \cap N(T)^\bot\cap \mbox{Dom}(T)$ 
and we observe that $h$ will stay in $N(T)^\bot$ for all times.
For $k=1$ the result is given by Theorem~\ref{theo:T}. 
We proceed by induction on $k$. 

First note that since the equation commutes with $x$-derivatives we have for 
the purely $x$-derivatives the analog of equation \eqref{eqx}, namely
  \begin{equation}\label{inequx:ho}
  \dt \|\partial_l ^0 h\|^2 _{L^2} \le - 2 \, \lambda \, \|\partial_l ^0 h - \Pi_l(\partial_l ^0 h) \|^2 _\Lambda.
  \end{equation}

For derivatives including some $v$-component ({\it i.e.}, $|j|\geq 1$), by means of {\bf H2'} with 
$\delta\leq (\nu_5 ^ \Lambda\nu_0 ^\Lambda)/(2\nu_1 ^\Lambda)$ we have the following estimate 
  \begin{multline}\label{inequv:ho}
  \dt \|\partial_l^j h\|^2 \leq 
  -2 \, \sum_{i \, | \, c_i(j) >0} \langle \partial_l^j h, \partial_{l+\delta_i}^{j-\delta_i}h\rangle_{L^2} 
  +\left(2 \, C(\delta) + 2 \, \nu^\Lambda _6 \right) \, \|h\|_{H^{k-1}}^2 
  - \nu^\Lambda _5 \, \|\partial_l^j h\|_\Lambda ^2.
  \end{multline}

For all $l$ with $|l|=k$ and $c_i(l)>0$, we consider the mixed term
  \begin{equation*}
  \dt \langle \partial_{l-\delta_i}^{\delta_i}h,\partial_l ^0 h\rangle = - \|\partial_l ^0 h\|^2 _{L^2} 
  +2 \, \langle L\left(\partial_l ^0 h-\Pi_l(\partial_l ^0 h)\right),\partial^{\delta_i}_{l-\delta_i} h\rangle. 
  \end{equation*}
By means of {\bf H1} we obtain 
  \begin{equation}\label{inequm:ho}
  \dt \langle \partial_{l-\delta_i}^{\delta_i}h,\partial_l ^0 h\rangle \leq-\|\partial_l ^0 h\|^2 _{L^2} 
  +C^L \, \eta \, \|\partial_l ^0h-\Pi_l(\partial_l ^0 h)\|^2 _\Lambda 
  +C^L \, \eta^{-1} \, \|\partial^{\delta_i}_{l-\delta_i}h\|^2_\Lambda. 
  \end{equation}

For $l$ with $|l|=k$ and $i$ such that $c_i(l)>0$, 
we define the following combination of derivatives of order $0$ and $1$ in $v$: 
  \begin{equation}
  \cal{Q}_{l,i}:=\alpha \, \|\partial_l ^0 h\|^2_{L^2} 
  + \beta \, \|\partial_{l-\delta_i}^{\delta_i}h\|^2 _{L^2} + 
  \gamma \, \langle \partial_{l-\delta_i}^{\delta_i}h,\partial_l h\rangle.
\end{equation}
By adjusting the constants $\alpha,\beta,\gamma>0$ and using Poincar\'e's inequality 
in the same way as in the proof of Theorem~\ref{theo:T} 
in Section~\ref{sec:main} it is straightforward to obtain 
  \begin{equation}\label{inequ:q}
  c\left(\|\partial_l ^0 h\|^2 _{L^2}+\|\partial_{l-\delta_i}^{\delta_i}h\|^2 _{L^2} \right)\leq \cal{Q}_{l,i}
  \leq C\left(\|\partial_l ^0 h\|^2 _{L^2} +\|\partial_{l-\delta_i}^{\delta_i}h\|^2 _{L^2}\right)
  \end{equation}
for some explicit constants $c,C>0$, and 
that the time derivative fulfills the following inequality
  \begin{equation}\label{inequ:dtq}
  \dt \cal{Q}_{l,i}\leq -K \, \left(\|\partial_l ^0 h\|^2 _\Lambda 
  +\|\partial_{l-\delta_i}^{\delta_i}h\|^2 _\Lambda \right)
  +C_0 \, \|h\|_{H^{k-1}}^2 
  \end{equation}
for some explicit constants $K,C_0>0$. 

Now we combine all the derivatives in the following way 
  \begin{equation}\label{Fk}
  \mathcal{F}_k(t):=\sum_{|l|+|j|=k, \, |j|\geq2} \left(\frac{\nu_0}{2}\right)^{-2|l|}\|\partial^j_l h\|^2
  +\frac{2}{K}\left(\frac{\nu_0}{2}\right)^{-2(k-1)}\sum_{|l|=k, \, i \, | \, c_i(l)>0} \cal{Q}_{l,i}
  \end{equation}
where $\nu_0=\nu_0^\Lambda/\nu_1^\Lambda >0$ is the constant such that (by assumption {\bf H1}) 
  \[ \|h\|_\Lambda \ge \nu_0 \, \|h\|_{L^2}. \]
By \eqref{inequ:q}, $\mathcal{\F}_k$ is equivalent to the square of the  
homogeneous Sobolev norm $\dot H^k$. 

To estimate the mixed terms in the right hand side of~\eqref{Fk} coming from~\eqref{inequv:ho}, 
we write 
  \begin{multline*}
  \sum_{|l|+|j|=k, |l|=s}2\left(\frac{\nu_0}{2}\right)^{-2s} \, 
  \langle \partial_l^j h, \partial_{l+\delta_i}^{j-\delta_i}h \rangle_{L^2} \leq\\
  \sum_{|l|+|j|=k, |l|=s}\left(\frac{\nu_0}{2}\right)^{-2s+1}\|\partial_l^j h\|^2
  +\sum_{|l|+|j|=k, |l|=s}\left(\frac{\nu_0}{2}\right)^{-2s-1}\|\partial_{l+\delta_i}^{j-\delta_i}h\|^2\leq\\
  \frac{1}{2}\sum_{|l|+|j|=k, |l|=s}
  \left(\frac{\nu_0}{2}\right)^{-2s}\|\partial_l^j h\|_\Lambda ^2+
  \frac{1}{2}\sum_{|l|+|j|=k,|l|=s}
  \left(\frac{\nu_0}{2}\right)^{-2(s+1)}\|\partial_{l+\delta_i}^{j-\delta_i} h\|_\Lambda ^2, 
  \end{multline*}
and we derive by combining \eqref{inequx:ho}, \eqref{inequv:ho} and \eqref{inequm:ho} that 
the time derivative of~\eqref{Fk} satisfies 
  \begin{multline*}
  \dt \mathcal{F}_k(t)\leq 
  C_+ \, \|h\|_{H^{k-1}}^2  
  -\sum_{|j|\geq2}\left(\frac{\nu_0}{2}\right)^{-2|l|}\|\partial_l^j h\|_{\Lambda} \\
  - \left(\frac{\nu_0}{2}\right)^{-2(k-1)}\left(\sum_{|j|=1}\|\partial_l^j  h\|^2 _\Lambda+
  2\sum_{|l|=k}\|\partial_l h \|^2 _\Lambda \right)\ .
  \end{multline*}
By using \eqref{inequ:q} we end up with
  \begin{equation*}
  \dt \mathcal{F}_k(t)\leq C_+ \, \|h\|_{H_{k-1}}^2 
  - K_- \, \left( \sum_{|j|+|l|=k}\|\partial^j _l h\|_\Lambda ^2 \right)
  \end{equation*}
for some explicit constants $C_+,K_->0$. Together with the induction 
assumption for $\mathcal{F}_1, \dots, \mathcal{F}_{k-1}$ this concludes 
the proof of the step $k$ by considering some 
combination of $\mathcal{F}_1, \dots, \mathcal{F}_k$. 
\end{proof}

\subsection{Weak external potential} 

If the particles are subject to an external force field, which is 
given as the gradient of a scalar potential $V$, the evolution equation on  
the distribution generalizes in the following way
  \begin{equation}\label{extpotequ}
  \partial_t f+v\cd\nab f-\nab V\cd\nabv f=Q(f,f).
  \end{equation}
We still consider for the spatial domain $\Omega = \tor$. 
We assume that $V=V(x)$ is $C^2$ and $e^{-V} \in L^1 _x$. 
For simplicity we restrict to 
collision operators for which the problem without potential admits 
a \emph{Maxwellian} equilibrium in this section and the next one. 
Moreover in this subsection we shall restrict further to collisional 
models admitting only mass conservation as a conservation law, {\it i.e.},  
we restrict to the relaxation model. 
The stationary solution is determined by mas conservation and given by
  \[ f_\infty = e^{-V} \, \frac{\rho_\infty}{\|e^{-V}\|_{L^1 _x}}  \, \cal{M} \]
where $\rho_\infty$ is the total mass of $f$ and $\cal{M}$ is the normalized Maxwellian
  \[ {\cal M}(v) = \frac{e^{-|v|^2/2}}{(2\pi)^{N/2}}. \]
We also assume that the external field is weak in the sense: 
  \begin{equation}\label{smallV}
  \| V \|_{C^2(\tor)} \le \var 
  \end{equation}
for some $\var$ depending on the collision operator.

Then we consider fluctuations around equilibrium of the form 
  \[ f = f_\infty + \sqrt{f_\infty} \, h \]
and we compute the following linearized equation on $h$: 
  \begin{equation}\label{eqlinV}
  \partial_t h +v\cd\nab h -\nab V\cd\nabv h = L(h) 
  \end{equation}
where $L$ is the linearized operator associated with $Q$ as 
before. We define the operator $T$ by 
  \[ T = L - v\cd\nab + \nab V\cd\nabv. \]
Let us assume that $L$ satisfies assumptions {\bf H1-H2-H3} and that 
the kernel of $L$ in $L^2 _v$ is given by $\mbox{Span} \{ \cal{M}^{1/2} \}$  
(this assumption is satisfied for the classical relaxation model). In this case the kernel 
of $T$ in $L^2 _{x,v}$ is trivially given by $\mbox{Span} \{ e^{-V/2} \cal{M}^{1/2} \}$. 
\medskip

We only sketch the proof of the following result:
 \begin{theorem}\label{theo:T:V}
 Under the assumptions {\bf H1-H2-H3} on $L$, 
 there is $\var_0>0$ such that for any $V \in C^1$ satisfying~\eqref{smallV} 
 with $\var \le \var_0$, the operator $T$ above satisfies the conclusion of Theorem~\ref{theo:T}. 
 If moreover conditions {\bf H1'-H2'-H3} hold for $L$,  and $V \in C^{k+1}$, then 
 the operator $T$ satisfies the conclusion of Theorem~\ref{theo:T:k}. 
 \end{theorem}

\begin{proof}
First let us recall that the $L^2$ norm of $h$ is decreasing using the non-positivity of $L$ 
and the skew symmetry of $- v\cd\nab + \nab V\cd\nabv$ in this Hilbert space. 

We will only show how to establish the bound on first-order derivatives. 
The generalization to higher-order is straightforward. 

Let us consider $h \in N(T)^\bot \cap \mbox{Dom}(T) \cap H^1 _{x,v}$. 
The time evolution for the $L^2$ norm of $h$ and the $L^2$ norm of its gradient in $v$ are unchanged. 
For the gradient in $x$, equation~\eqref{eqx} is replaced by
  \begin{equation*}
  \dt \|\nab h\|_{L^2}^2\leq - 2 \, \lambda \, \|\nab h-\Pi_l(\nab h)\|_\Lambda ^2+ 
   2 \, \var \, \| \nabla_v h\|_{L^2} ^2 \|\nabla_x h\|_{L^2} ^2.
  \end{equation*}

Finally for the mixed term we have to replace \eqref{eqxv} by
 \begin{multline*}
  \dt \langle \nabla_x h, \nabla_v h \rangle \le - \left\|\nabla_x h\right\|^2 _{L^2} 
             + C^L \, \eta \, \left\| \nabla_x h - \Pi_l(\nabla_x h) \right\|^2 _\Lambda \\
              + C^L \, \eta^{-1} \, \left\| \nabla_v h \right\|^2 _\Lambda+\var\|\nabla_v h\|^2_{L^2}.
 \end{multline*} 

Hence following the same arguments as in Section~\ref{sec:main}, 
we find the following differential inequality on the quadratic form $\cal{F}$: 
 \[
 \frac{{\rm d} {\cal F}(t)}{{\rm d}t} \le 
             - K ( \|h\|_\Lambda ^2 + \|\nabla_{x,v} h \|_\Lambda ^2 ) 
           + a \, \var \left\{\, \| \nabla_v h\|_{L^2}  \|\nabla_x h\|_{L^2}+\|\nabla_v h\|^2_{L^2} \right\} 
 \]
for some explicit constant $a,K>0$. Therefore it concludes the proof for $\var>0$ small enough. 
%
\end{proof}

\Remarks 

1. It may be possible that for a spatial domain 
$\Omega = \R^N$ a modified version of this strategy could be applied, 
assuming additionally that $V$ satisfies a log-Sobolev inequality 
on $\Omega$. 
\smallskip 

2. This subsection about weak external fields illustrates the fact that our method 
is robust, since it is based on {\em a priori} estimates, which remain true up to a perturbation. 
\medskip
%


\subsection{Self-consistent potential}

Let us consider a collisional kinetic model for 
particles which interact through collisions and also 
through a self-consistent potential. In this subsection we exclude the semi-classical 
relaxation collision operators. For the 
potential we consider the physically most common case of Poisson interaction. 
More precisely 
  \begin{equation}\label{eqself}
  \left\{
    \begin{array}{l}
    \partial_t f +v\cd\nab f- \epsilon \, \nab V \cd\nabv f =Q(f,f) \vspace{0.2cm} \\
    \displaystyle 
    \Delta_x V= \rho - \rho_0, \quad \int_{\ens{T}^N} V \, \ud x =0,   
    \end{array}
  \right.
  \end{equation}
where the coupling is {\it via} the density $\rho(t,x)=\int f(t,x,v)\, \ud v$.
The equation models particles interacting by binary collision or by scattering a 
background medium in thermal equilibrium, and at the same time 
interacting {\it via} electrostatic forces in the case $\epsilon =-1$ or by gravitational attraction 
in the case $\epsilon = +1$. 
Existence of global classical solutions in the large to the 
Vlasov-Poisson system in the torus has been proven in \cite{BaRe}. 
For the system including various collision operators, 
the existence of classical solutions has been established very recently 
in the articles discussed in the introduction. 

We consider the previous equation in the torus $x \in \tor$ for $\epsilon=-1$ 
(electrostatic interaction). It admits a global Maxwellian equilibrium $f_\infty$, 
determined by the conservation laws (for a more 
general discussion of the possible equilibria, we refer to~\cite{DeDo91} for 
instance). 

Then we consider the linearization $f=f_\infty + \sqrt{f_\infty} h$ around this equilibrium. 
Discarding the bilinear terms yields 
  \begin{equation*}
  \left\{
    \begin{array}{l} \displaystyle 
    \partial_t h +v\cd\nab h - (v \cdot \nab V) \, f_\infty ^{1/2} = L(h)  \vspace{0.2cm} \\ \displaystyle 
    \Delta_x V= \left( \int_{\R^N} h \, f_\infty ^{1/2} \, \ud v \right)  
    \end{array}
  \right.
  \end{equation*}
where $L$ is the linearized collision operator associated with $Q$. 

Now let us denote by $T_p$ the operator on $L^2$ defined by  
  \begin{equation}\label{def:T:Poisson}
  T_p= L(h) - v \cdot \nabla_xh + (v \cd\nabla_x V(h)) \, f_\infty ^{1/2}.
  \end{equation}

Then (defining $\Pi_g$ as before) we have the following theorem
 \begin{theorem}\label{theo:T:Poisson}
 Let $L$ satisfy the assumptions {\bf H1-H2-H3}. 
 Then the conclusion of Theorem~\ref{theo:T} still holds true for the operator $T_p$ defined in~\eqref{def:T:Poisson}. 
 If moreover $L$ satisfies assumptions {\bf H1'-H2'-H3}, then the operator 
 $T_p$ obeys the conclusion of Theorem~\ref{theo:T:k}. 
 \end{theorem}

\begin{proof}[Proof of Theorem~\ref{theo:T:Poisson}]
The proof is almost exactly the same as to the one of Theorem~\ref{theo:T}. Therefore 
we shall only indicate the differences in the estimates. Essentially  
the norm has to be modified in order to take into account the interaction 
energy. For the $L^2$ norm, one has by integration by parts and using that $L$ is mass conserving
  \[ \frac{\ud}{\ud t} \left( \| h \|_{L^2} ^2 + \| \nabla_x V \|_{L^2 _x} ^2 \right) 
      \le  - 2 \, \lambda \, \| h - \Pi_l (h) \|_\Lambda ^2. \]

Similarly one has on the gradient in $x$ 
  \[ \frac{\ud}{\ud t} \left( \| \nabla_x h \|_{L^2} ^2 + \sum_{1 \le i,j \le N} 
      \left\| \frac{\partial^2 V}{\partial x_i \partial x_j} \right\|_{L^2 _x} ^2 \right) 
      \le - 2 \, \lambda \, \| \nabla_x h - \Pi_l (\nabla_x h) \|_\Lambda ^2. \]

For the time evolution of the gradient in $v$ one has the additional term 
  \[ -2 \, \int \left( \nabla_x V \cdot \nabla_v f_\infty ^{1/2}\right) \, h \]
and for the mixed term, one has the additional terms 
   \[ - 2 \, \int \Delta_x V \, h \, f_\infty ^{1/2} 
       + C \,  \int \left( \sum_{i,j} v _i v _j \frac{\partial^2 V}{\partial x _i \partial x_j} \right) 
         \, h \, f_\infty ^{1/2}.  \]
Since it is straightforward that all these additional terms as well as 
the gradient of $V$ can be controlled by the $L^2$ 
norm of $h$, the end of the proof is straightforward as in Theorem~\ref{theo:T}. 
\end{proof}

\Remarks

 1. We need the interaction to be repulsive in order to have 
the right sign for the interaction energy. For gravitational self-interaction 
potentials, our method would work as well under an additional assumption of 
smallness of the potential, {\it i.e.}, assuming $\epsilon >0$ and $\epsilon \le \epsilon_0$ 
with $\epsilon_0$ depending on the collision operator in~\eqref{eqself}   
(the proof is similar to the case of a weak external potential).  
\smallskip

2. As noticed in~\cite{Guo:PB:02}, the non-linear term arising from a 
self-interaction Poisson potential can be controlled, in the energy estimates 
for $H^k$ with $k$ such that $E(k) >N/2$, by 
  \[ C \, \|h\| _{H^k} \, \|h \langle v \rangle^{1/2} \|_{H^k} ^2, \]
using Sobolev embeddings and the straightforward elliptic estimate 
in the torus: 
  \[ \forall \, l \ge 1, \qquad \|\nabla_x V\|_{H^l} \le C_l \, \|f\|_{H^{l-1}}. \]  
Therefore one can extend the construction of
smooth solutions near equilibrium in Theorem~\ref{theo:NL} to the case 
when a self-consistent Poisson potential is added, as long as the coercivity 
norm $\Lambda$ of the linearized problem is stronger than the norm 
$\|h \langle v \rangle^{1/2} \|_{L^2}$. This is true, for instance, for 
the Boltzmann equation for hard spheres (considered in~\cite{Guo:PB:02}), and also for the Landau equation 
with $\gamma \ge -1$ (see Subsection~\ref{ss:La}). 
\medskip

\section{Application to full non-linear models near equilibrium}\label{sec:NL}
\setcounter{equation}{0}

In this section we prove existence and uniqueness of 
smooth global solutions near equilibrium thanks to the coercivity estimates on
the linearized models. This yields also explicit exponential rate of convergence to equilibrium. 
Obviously there is nothing else to prove for linear models (such as the classical relaxation or 
linear Fokker-Planck equation). 
Therefore let us assume that the collision operator is bilinear and 
let us denote the remaining term in the linearization process 
$f= f_\infty + f_\infty ^{1/2} h$: 
  \[ \Gamma(h,h) = f_\infty ^{-1/2} \, Q(f_\infty h , f_\infty h). \]

In this section we shall consider a linearized collision operator $L$ 
satisfying assumptions {\bf H1'-H2'-H3} 
(we make the additional assumptions {\bf H1'-H2'} in order to get 
coercivity estimates in higher-order Sobolev spaces). Moreover we shall 
assume on the bilinear form $\Gamma$:
 \begin{itemize}
 \item[{\bf H4.}] There is $k_{0} \in \N$ and $C_{\Gamma} >0$ such that for $k \ge k_{0}$,  
   \[ \| \Gamma(h,h) \|_{H^k} \le C_\Gamma \, \| h \|_{H^k} \, 
           \left( \sum_{|j|+|l| \le k} \|\partial^j _l h\|_\Lambda ^2 \right)^{1/2}. \]
 \end{itemize}

Then we have 
 \begin{theorem}\label{theo:pert} 
 Let $Q$ a (bilinear) collision operator such that 
 \begin{itemize}
 \item[(i)] equation~\eqref{eqgenNL} admits an equilibrium 
 $0 \le f_\infty \in L^1(\tor \times \R^N)$; 
 \item[(ii)] the linearized collision operator $L=L(h)$ around $f_\infty$ with the 
 scaling $f=f_\infty + f_\infty ^{1/2} h$ satisfies {\bf H1'-H2'-H3};
 \item[(iii)] the bilinear remaining term $\Gamma=\Gamma(h,h)$ in the linearization 
 satisfies {\bf H4}. 
 \end{itemize}
 
 Then for any $k \ge k_{0}$ (where $k_{0}$ is defined in {\bf H4}), 
 there is $\var_0>0$ such that for any distribution $0 \le f_0 \in L^1$ with
   \[ \| (f_0 - f_\infty) \, f_\infty ^{-1/2} \|_{H^k} \le \var_0 \]
 there exists a unique global smooth solution 
 $0 \le f=f(t,x,v)$ to equation~\eqref{eqgenNL}, which satisfies 
   \[ \| (f_t - f_\infty) \, f_\infty ^{-1/2} \|_{H^k} \le C_0 \, \var_0 \, e^{-\tau t} \]
 for some explicit constant $C_0,\var_0, \tau >0$, depending only on the constants 
 appearing in {\bf H1'-H2'-H3-H4}.
 \end{theorem}

\begin{proof}[Proof of Theorem~\ref{theo:pert}]
We explain the proof by {\em a priori} arguments, on a given 
smooth solution. The construction of positive solutions thanks to the estimates 
above is based on, by now standard, fixed point arguments 
(we refer the reader to \cite{Guo:PB:02,Guo:LE:02} for instance).

The function $h = (f - f_\infty) \, f_\infty ^{-1/2}$ satisfies $\Pi_{g}(h)=0$ and it solves 
  \[ \partial_t h = T(h) + \Gamma(h,h). \]
Then we estimate the time evolution of the $\mathcal{H}^k$ norm, defined in Theorem~\ref{theo:T:k}: 
  \[ \dt \|h\|_{\mathcal{H}^k} ^2 = 2 \langle Th, h \rangle_{\mathcal{H}^k} + 2 \langle 
     \Gamma(h,h), h \rangle_{\mathcal{H}^k}. \]
We deduce that 
  \[ \dt \|h\|_{\mathcal{H}^k} ^2 \le - C_{T} \, \left( \sum_{|j|+|l|\le k} \|\partial^j _l h\|_\Lambda ^2 \right) 
     + C_{\var} \, \| \Gamma(h,h) \|_{H^k} ^2 + \var \, \| h \|_{H^k} ^2. \]

Hence, by taking $\var$ small enough,
  \[ \dt \|h\|_{\mathcal{H}^k} ^2 \le - \frac{C_{T}}2 \, 
            \left( \sum_{|j|+|l|\le k} \|\partial^j _l h\|_\Lambda ^2 \right) 
     + C_\Gamma ' \, \| h \|_{H^k} ^2 \, \left( \sum_{|j|+|l| \le k} \|\partial^j _l h\|_\Lambda ^2 \right). \]

This concludes the proof by maximum principle since the $\Lambda$ norm 
controls the $L^2$ norm. 
\end{proof}

\section{Proof of the general assumptions for physical models}\label{sec:phys}
\setcounter{equation}{0}

\subsection{Linear relaxation}

We consider the linear relaxation equation in the torus
 \begin{equation}\label{eqf}
 \partial_t f + v \cdot \nabla_x f = \frac{1}{\kappa} \left[ 
    \left( \int_{\R^N}  f(t,x,v_*) \, \ud v_* \right) {\cal M}(v) -f \right],
 \end{equation}
for $x \in \ens{T}^N$ and $v \in \R^N$ ($N \ge 1$). 
Here $\kappa >0$ denotes the {\em Knudsen number} and 
${\cal M}$ denotes the normalized Maxwellian:
 \[ {\cal M}(v) = \frac{e^{-|v|^2/2}}{(2\pi)^{N/2}} \]
with mass $1$, momentum $0$ and temperature $1$. 
This equation preserves the total mass of the distribution 
  \[ \forall \, t \ge 0, \ \ \ \int_{\ens{T}^N \times \R^N} 
       f(t,x,v) \, \ud x \, \ud v = \int_{\ens{T}^N \times \R^N} f_0(x,v) \,\ud x \, \ud v \]
but admits no other conservation law. 
For a given initial datum $f_0 \ge 0$, it admits a unique 
global equilibrium $f_\infty=\rho_\infty \, \cal{M}$, where $\rho_\infty$ 
is the total mass of $f_0$, defined by 
  \[ \rho_\infty = \int_{\ens{T}^N \times \R^N} f_0(x,v) \, \ud x \, \ud v. \] 
Finally let us add that the Cauchy theory is straightforward for \eqref{eqf} 
since it is linear (see~\cite{CCG03} for instance for more details on this equation). 

We rescale the equation as
  \[ f = f_\infty + \sqrt{f_\infty} h = \rho \, {\cal M} + \rho^{1/2} M \, h, \]
where $M := \sqrt{{\cal M}}$. The equation for $h$ reads 
  \begin{equation*}
  \partial_t h + v \cdot \nabla_x h = \frac{1}{\kappa} \, \left[ 
      \left( \int_{\R^N} h' \, M' \, \ud v' \right) M - h \right] =: L(h) 
  \end{equation*}
where we have used the classical notation $h'=h(v')$. 
We split the operator $L$ into 
  \[ L = K - \Lambda, \quad K(h) = \frac{1}{\kappa} \, \left( \int_{\R^N} h' \, M' \, \ud v' \right) M, 
      \quad \Lambda(h) = \kappa^{-1} \, h. \] 
Therefore $L$ satisfies {\bf H1} taking $\|\cdot\|_\Lambda = \|\cdot \|_{L^2 _{x,v}}$, 
and assumption {\bf H2} follows straightforwardly (with $C(\delta)=0$) from 
   \[ \nabla_v K(h) = \frac{1}{\kappa} \, \left( \int_{\R^N} h' \, M' \, \ud v' \right) \nabla_v M. \]
Observe also that the strengthened assumptions {\bf H1'} and {\bf H2'}, 
that are necessary to ensure decay in higher-order Sobolev norms, are satisfied straightforwardly.

The operator $L$ is local in $x$ and $t$. When $x$ is fixed, it is well-defined and bounded 
on $L^2 _v$, and it is self-adjoint non-positive on this space. More precisely 
its Dirichlet form is given by  
  \[ \langle L(h), h \rangle_{L^2 _v} = - \frac{1}{2 \kappa} \, 
      \int_{\R^N \times \R^N} \left( \frac{h}{M} -\frac{h'}{M'} \right)^2 
        \, {\cal M} \, {\cal M'} \, \ud v \, \ud v'. \] 
Therefore its kernel is $N(L) = \mbox{Span}\left\{M \right\}$, and 
we define $\Pi_l$ the (orthogonal) projection on this space in $L^2 _v$: 
  \[ \Pi_l (h) = \left( \int_{\R^N} h' \, M' \, \ud v' \right) M = \kappa \, K(h). \]
Then it is straightforward that $L$ has a spectral gap $\lambda_L = \kappa^{-1}$, since 
  \[ \int_{\R^N} L(h) \, h \, \ud v = - \frac{1}{\kappa} \, \|h - \Pi_l(h)\|_{L^2} ^2. \]
Hence assumption {\bf H3} is satisfied. 

\subsection{Semi-classical relaxation}

Let us now modify the previous relaxation equation in order to 
take into account quantum effects associated with the Pauli's exclusion principle. 
We thus consider the following semi-classical model (weakly non-linear) 
for charged particles (here $\epsilon \in \{-1,0,1\}$): 
  \begin{equation}\label{eqfsc}
  \partial_t f + v \cdot \nabla_x f 
   = \frac{1}{\kappa} \, 
   \int_{\R^N} \big[ {\cal M}(1-\epsilon f) f' - {\cal M}' (1-\epsilon f') f \big] \, \ud v' 
   := Q_\epsilon (f,f),
  \end{equation}
for $x \in \ens{T}^N$ and $v \in \R^N$, 
where $\mathcal{M}$ is the normalized Maxwellian as before, and 
we have used the shorthand $f'=f(v')$. Regardless of the choice of 
$\epsilon$ this scattering operator is mass preserving but 
admits no other conservation law. The case $\epsilon=0$ 
is the standard linear relaxation model of the previous subsection. 

In the case $\epsilon=1$ this operator is probably the simplest model 
describing a gas of fermions relaxing towards the thermodynamic equilibrium 
for a perfect fermigas, that is the {\em Fermi-Dirac distribution}. 
The operator $Q_{+1}$ describes 
the interaction of the fermions with a background medium at rest with  
constant temperature. The factors $(1-f)$ correspond to correlation of 
particles before and after collision due to Pauli's exclusion principle. 
Note that this modification of the standard relaxation mechanism is (at 
least for the usual range of temperatures and densities) necessary only in 
particular situations, such as the one of a gas of electrons. 
Those are, due to their small mass, most likely 
to satisfy Sommerfeld's degeneracy condition (see~\cite{CC90}). This 
equation can be seen as the scattering counterpart of the full fermion 
Boltzmann equation studied for example in~\cite{CC90,Dol94}. 
For a more detailed introduction to models describing scattering as well 
as binary collisions for fermions see \cite{MRS90}. 
For this model, a Cauchy theory can be obtained 
using maximum principle arguments to treat the (weak) nonlinearity 
(see~\cite{NS04}) assuming some bounds on the initial datum.
The long-time behaviour of solutions to this equation has been studied 
by the EEP-method in \cite{NS04} (leading to polynomial rates of convergence 
to equilibrium), however the necessary uniform regularity 
bounds on the solution were assumed. 

In the case $\epsilon=-1$, $Q_{-1}$ describes the interaction of 
bosons with background medium at rest with some constant temperature, 
and is probably the simplest model 
describing a gas of bosons relaxing towards the thermodynamic equilibrium 
for a perfect bosongas, that is the {\em Bose-Einstein distribution}. 
In this setting the existence of a boson 
in a velocity and space interval will increase the chance of another 
boson being scattered to this interval (see~\cite{CC90}). This mechanism 
leads to Bose-Einstein condensation for low temperatures and large densities. 
However since we linearize around the \emph{regular} equilibrium we cannot describe this phenomenon but 
only the situation farther away from the critical mass for the phase 
transition. In the spatially homogeneous setting very precise asymptotics 
including optimal rates for the convergence above as well as below the critical 
mass have been given in \cite{EMV}.
The authors study a model for Compton scattering of photons against electrons. In the non condensate case their model - in an appropriate scaling - is analogous to the one we study here because in this case their cross-section is bounded away from $0$. The authors prove exponential convergence in the non-condensate case, 
which is consistent with the result we obtain in the $x$-dependent situation.
Moreover they derive an optimal rate for the convergence in the condensate case which is only polynomial. Thus is seems unavoidable to impose a bound on the initial mass to retain exponential convergence.
The existence of solutions to a more elaborate collisional model (in the spatially 
homogeneous case) has been shown in \cite{Lu00} (see also~\cite{Lu04}). In the 
same work the weak convergence to regular equilibrium states (similar to 
the ones that our simplified model admits) has been shown to hold true above 
a specific temperature, larger than the critical one, and for which an explicit bound is given.

The equation preserves the total mass of the distribution 
  \[ \forall \, t \ge 0, \ \ \ \int_{\ens{T}^N \times \R^N} 
       f(t,x,v) \, \ud x \, \ud v = \int_{\ens{T}^N \times \R^N} f_0(x,v) \, \ud x \, \ud v = \rho\]
and admits (recalling in the boson case the bound imposed on the initial 
datum in Theorem~\ref{theo:NLSC}) a unique equilibrium 
  \[ f_\infty = \frac{\kappa_\infty {\cal M}}{1 + \epsilon \kappa_\infty {\cal M}}, \]
where $\kappa_\infty$ is determined by the mass $\rho$ of $f_0$. 

We linearize the equation for a general scaling $f = f_\infty + m \, h$ 
($m$ is a given positive function only depending on $v$). Discarding the bilinear term, 
the equation for $h$ reads after straightforward computations
  \begin{equation*}
  \partial_t h + v \cdot \nabla_x h = \frac{1}{\kappa} \, 
  \int_{\R^N} \left[ h' \, \frac{m' (1+\epsilon \kappa_\infty {\cal M}')}{m (1+\epsilon 
    \kappa_\infty {\cal M})} {\cal M} 
   - h \, \frac{1+\epsilon \kappa_\infty {\cal M}}{1+ \epsilon \kappa_\infty {\cal M}'} {\cal M}' \right] 
\, \ud v' =: L_m(h). 
  \end{equation*}
We make the choice $m= (1+\epsilon \kappa_\infty {\cal M})^{-1} \, M\sqrt{\kappa_{\infty}}$, 
where $M:=\sqrt{{\cal M}}$. It yields $L = K - \Lambda$ with 
  \[ \quad K(h) = \frac{1}{\kappa} \, \left( \int_{\R^N} h' \, M' \, \ud v' \right) M, \] 
and $\Lambda$ is the multiplicative operator by $\nu$ with 
  \[ \nu(v)= \frac{1}{\kappa} \, 
      (1+\epsilon \kappa_\infty {\cal M}) \, 
     \left( \int_{\R^N} \frac{{\cal M}'}{1+ \epsilon \kappa_\infty {\cal M}'} \, \ud v' \right) 
      = \frac{\rho}{\kappa \kappa_{\infty}} \, (1+\epsilon \kappa_\infty {\cal M}). \]
Therefore $L$ satisfies {\bf H1} taking $\|\cdot\|_\Lambda = \|\cdot\|_{L^2 _{x,v}}$, and
assumption {\bf H2} follows straightforwardly (with $C(\delta)=0$) from 
  \[ \nabla_v K(h) = \frac{1}{\kappa} \, \left( \int_{\R^N} h' \, M' \, \ud v' \right) \nabla_v M.  \]
Again the strengthened assumptions {\bf H1'} and {\bf H2'} are also satisfied straightforwardly.

The resulting operator $L$ is local in $x$ and $t$. When $x$ is fixed, it is well-defined and bounded 
on $L^2 _v$, and it is self-adjoint non-positive on this space. 
More precisely its Dirichlet form is given by 
  \begin{multline*}
    \langle L(h), h \rangle_{L^2} =  -\frac{1}{2\kappa} \, 
    \int_{\R^N \times \R^N} \left( \frac{h(1+\epsilon \kappa_\infty {\cal M})}{M} 
      -\frac{h'(1+\epsilon \kappa_\infty {\cal M}')}{M'} \right)^2 \times \\ 
       \mbox{ } \hspace{7cm}  (1+\epsilon \kappa_\infty {\cal M})^{-1} \, 
     (1+\epsilon \kappa_\infty {\cal M}')^{-1} \, {\cal M} \, {\cal M'} \, \ud v \, \ud v' \\
  = -\frac{1}{2\kappa} \, 
    \int_{\R^N \times \R^N} \left( \frac{h M}{f_\infty} 
      -\frac{h' M'}{f_\infty '} \right)^2 \, f_\infty \, f_\infty ' \, \ud v \, \ud v'.
  \end{multline*}
Therefore its kernel is 
  \[ N(L) = \mbox{Span}\left\{\frac{f_\infty}{M} \right\}. \]
We define $\Pi_l$ the (orthogonal) projection on this space in $L^2 _v$:
  \[ \Pi_l (h) = \left( \int_{\R^N} h' \, \frac{f_\infty '}{M'} \, \ud v' \right) 
          \frac{f_\infty}{M}. \]

First let us explain how to show assumption {\bf H3} by non-constructive approach. 
As $\nu \ge \underline{\nu}$ with $\underline{\nu}= \inf_{\mathbb{R}^N}\nu>0$, 
we deduce that the (bounded) multiplicative operator $\Lambda$ on $L^2_v$ has its spectrum 
included in $(-\infty,-\underline{\nu}]$. Then the operator $K$ is straightforwardly compact 
on $L^2 _v$ and thus by Weyl's theorem about compact perturbation of the essential 
spectrum of a self-adjoint operator in a Hilbert space, we deduce that the essential 
spectrum of $L$ is included in $(-\infty,-\underline{\nu}]$. The remaining discrete spectrum 
lies in $\R_-$ because of the sign of the Dirichlet form, and as $0$ is an isolated 
eigenvalue, we deduce that there is $\lambda_0 >0$ such that the non-zero part of 
the spectrum of $L$ lies in $(-\infty,-\lambda_0]$. This immediately shows {\bf H3} 
with $\lambda_L = \lambda_0/\overline{\nu}$ where $\overline{\nu} = \sup_{\R^N} \nu < +\infty$. 

It is interesting to note that if the mass in the boson case matches 
the critical one then $\nu$ is still nonnegative but it becomes zero for $v=0$, 
which of course is the place where condensation happens. 
In this case our method breaks down because we lack the spectral gap (the essential 
spectrum reaches $0$ and one only expects some polynomial rates of convergence to equilibrium, 
which is consistent with~\cite{EMV}).
\smallskip

Second we restrict to the fermionic case ($\epsilon=1$) and
we explain how to estimate explicitly $\lambda_L$. 
Let us consider some function $h$ orthogonal to $M^{-1} f_{\infty}$. Then 
  \[ \int_{\R^N} h M \, \ud v = \int_{\R^N} h \left( M - \frac{f_{\infty}}{\kappa_{\infty}M} \right) \, \ud v 
                    =  \, \int_{\R^N} h 
                  \frac{\kappa_{\infty}M^3}{1+  \kappa_{\infty}{\cal M}} \, \ud v. \] 
Hence we deduce that 
  \begin{multline*} 
  \left( \int_{\R^N} h M \, \ud v \right)^2 \le \left( \int_{\R^N} h^2 \nu \, \ud v \right) 
      \left( \int_{\R^N} \nu^{-1} \, 
        \frac{\kappa_{\infty} ^2 M^6}{(1+ \kappa_{\infty}{\cal M})^2} \, \ud v \right) \\
        = \frac{\kappa}{\rho} \, \left( \int_{\R^N} h^2 \nu \, \ud v \right) 
      \left( \int_{\R^N} f_{\infty} ^3 \, \ud v \right).
  \end{multline*}
From the exact formula for $f_\infty$ it is straightforward that 
$f_\infty ^3 < f_\infty$ and thus
  \[ \frac{1}{\rho} \, \left( \int_{\R^N} f_{\infty} ^3 \, \ud v \right) < 1. \] 
We get 
	  \begin{multline*}
          \langle L(h), h \rangle_{L^2} = 
          \left[ \frac1\kappa \,  \left( \int_{\R^N} h M \, \ud v \right)^2 - 
         \left( \int_{\R^N} h^2 \nu \, \ud v \right) \right] \\ 
          \le -\left[ 1 - \frac{1}{\rho} \, \left( \int_{\R^N} f_{\infty} ^3 \, \ud v \right) \right] 
                 \left( \int_{\R^N} h^2 \nu \, \ud v \right).  
          \end{multline*}
Since in the fermionic case $\nu \ge (\kappa \kappa_{\infty})^{-1}$, we deduce that $L$ has 
a spectral gap $\lambda_0$, with the explicit estimate:
 	\[ \lambda_{L} \ge \frac{1}{\kappa \kappa_{\infty}} \, 
         \left[ 1 - \frac{1}{\rho} \, \left( \int_{\R^N} f_{\infty} ^3 \, \ud v \right) \right]. \]
Hence we deduce an explicit estimate on $\lambda_L$ since $\lambda_L = \lambda_0/\overline{\nu}$ 
where $\overline{\nu} = \sup_{\R^N} \nu= (1 + \kappa_\infty \cal{M}(0) )/(\kappa \kappa_\infty)$ 
is explicit. 

Now we want to establish the bound {\bf H4} on the bilinear part. It is given by
\begin{equation*}
\Gamma(h,h)=
\frac{\epsilon \kappa_\infty ^{1/2}}{\kappa} \, h \, 
\left( \int_{\R^N} h' \, M' \, \frac{(\cal{M'} - \cal{M})}{1+\kappa_\infty \cal{M'}} \, \ud v' \right).
\end{equation*}
Therefore {\bf H4} is immediately obtained by using Leibniz rule on higher-order 
derivatives, the trivial bound $L^2 \times L^2 \to L^2$ on $\Gamma$, and 
Sobolev embeddings (which requires that $E(k_0/2) > N/2$ where $E$ denotes 
the entire part of a real number). 

\Remark The scaling that we used to linearize the collision operator is not exactly 
the same as in Theorem~\ref{theo:pert} since we 
choose $m=f_\infty^{1/2}\left(1+\epsilon\kappa_\infty\mathcal{M}\right)^{-1/2}$. 
However it is easy to see (following exactly the same proof) that the statement 
of Theorem~\ref{theo:pert} remains true also with the scaling $f=f_\infty+m \, h$ 
when the factors $f_\infty^{-1/2}$ are replaced by the some factors $m^{-1}$ 
with the same decay at large velocities. This leads to the statement of Theorem~\ref{theo:NLSC}.
\medskip
 
%
%
%
%

\subsection{The linear Fokker-Planck equation}

We consider the linear Fokker-Planck equation in the torus
 \begin{equation}\label{eqFP}
 \partial_t f + v \cdot \nabla_x f = \nabla_v \cdot \left( \nabla_v f + fv \right),
 \end{equation}
for $x \in \ens{T}^N$ and $v \in \R^N$ ($N \ge 1$). 

This equation preserves the total mass of the distribution 
  \[ \forall \, t \ge 0, \ \ \ \int_{\ens{T}^N \times \R^N} 
       f(t,x,v) \, \ud x \, \ud v = \int_{\ens{T}^N \times \R^N} f_0(x,v) \,\ud x \, \ud v \]
but admits no other conservation law. 
For a given initial datum $f_0 \ge 0$, it admits a unique 
global equilibrium $f_\infty=\rho_\infty \, \cal{M}$, where $\rho_\infty$ 
is the total mass of $f_0$, defined by 
  \[ \rho_\infty = \int_{\ens{T}^N \times \R^N} f_0(x,v) \, \ud x \, \ud v, \] 
and ${\cal M}$ is the normalized Maxwellian distribution 
 \[ {\cal M}(v) = \frac{e^{-|v|^2/2}}{(2\pi)^{N/2}} \]
with mass $1$, mean $0$ and temperature $1$. 

We study fluctuations around the equilibrium in the form
  \[ f = f_\infty + \sqrt{f_\infty} h = \rho_\infty \, {\cal M} + \rho_\infty ^{1/2} M \, h, \]
where $M := \sqrt{{\cal M}}$. The equation for $h$ reads 
  \begin{equation}\label{eqFPh}
  \partial_t h + v \cdot \nabla_x h = \Delta_v h + \left( \frac{N}2 - \frac{|v|^2}4 \right) \, h =: L(h).  
  \end{equation}
We split the operator $L$ into 
  \[ L = K - \Lambda, \quad K =0, 
      \quad \Lambda = L \]
($K=0$ is typical for a purely diffusive collisional model). 
Assumption {\bf H2} is obviously fulfilled with $C(\delta)=0$. 

Let us prove that $L$ satisfies {\bf H1} taking 
  \[ \| h \|_\Lambda = \left( \| v \, h \|_{L^2} ^2 + \|\nabla_v h \|_{L^2} ^2 \right)^{1/2}. \] 
Indeed one can check easily that this norm is stronger than the $L^2$ norm: 
  \[ \| h \|_\Lambda ^2 \ge 2 \, \| v \, h \|_{L^2} \, \|\nabla_v h \|_{L^2} 
                \ge - 2 \, \int_{\R^N} (v\,h) \cdot \nabla_v h \, \ud x \, \ud v 
                = N \, \| h\|_{L^2} ^2. \]
Straightforward computations yields 
  \begin{equation}\label{eq:FP1}
  \langle \Lambda  (h), h \rangle_{L^2} \ge C_1 \, \| h \|_\Lambda ^2 - C_2 \, \|h\|_{L^2} ^2 
  \end{equation}
for explicit constants $C_1, C_2 >0$. The operator $L$ is local in 
$x$ and $t$. When $x$ is fixed, it is well-defined and bounded 
on $L^2 _v$, and it is self-adjoint non-positive on this space. More precisely 
its Dirichlet form is given by  
  \begin{equation*}
  \langle L(h), h \rangle_{L^2 _v} = - 
         \int_{\R^N} \left\| \nabla_v h + \frac{v}2 \, h \right\|^2 \, \ud v. 
  \end{equation*} 
Therefore its kernel is $N(L) = \mbox{Span}\left\{M \right\}$, and 
we define $\Pi_l$ the (orthogonal) projection on this space in $L^2 _v$: 
  \[ \Pi_l (h) = \left( \int_{\R^N} h' \, M' \, \ud v' \right) M. \]
Classical computations based on Poincar\'e's inequality with measure 
$\cal{M}$ show that 
  \begin{equation}\label{eq:FP2}
  \int_{\R^N} \left\| \nabla_v h + \frac{v}2 \, h \right\|^2 \, \ud v  \ge 2 \, \| h \|_{L^2} ^2. 
  \end{equation} 
Then combining~(\ref{eq:FP1},\ref{eq:FP2}) yields 
  \[ \langle L(h), h \rangle_{L^2} \le - \lambda \, \| h \|_\Lambda ^2 \]
for some explicit constant $\lambda >0$.  

Finally we have 
  \begin{multline*}
  \langle \nabla_v L(h), \nabla_v h \rangle_{L^2} = 
      \langle L( \nabla_v h), \nabla_v h \rangle_{L^2} - \langle (v/2) h, \nabla_v h \rangle_{L^2} \\
      = \langle L( \nabla_v h), \nabla_v h \rangle_{L^2} + \frac{N}2 \, \|h\|_{L^2} 
      \le - \Lambda_L \, \| h \|_\Lambda ^2 + \frac{N}2 \, \|h\|_{L^2}. 
  \end{multline*}
The two last inequalities conclude the proof of {\bf H1} and {\bf H3}. 


\subsection{The Boltzmann equation}\label{linbol}

Let us consider the Boltzmann equation (here $N \ge 2$)
  \begin{equation}\label{eq:Bol}
  \partial_{t} f + v \cdot \nabla_{x} f = Q(f,f), \quad t \ge 0,\ x \in \tor, \ v \in \R^N
  \end{equation} 
with a collision operator (local in $t,x$)
  \[ Q(f,f) = \int_{\R^N \times \ens{S}^{N-1}} B(|v-v_{*}|, \cos \theta) \, 
                    \left( f' f'_{*} - f f_{*} \right) \, \ud v_{*} \, \ud\sigma. \]
We adopt the notations $f' = f(v')$, $f_{*} = f(v_{*})$ and $f' _{*} = f(v' _{*})$, where
 \begin{equation*}\label{eq:rel:vit}
 v' = (v+v_*)/2 + (|v-v_*|/2) \, \sigma, \qquad
 v'_* = (v+v_*)/2 - (|v-v_*|/2) \, \sigma
 \end{equation*}
stand for the pre-collisional velocities of particles which after
collision have velocities $v$ and $v_*$. Moreover $\theta\in
[0,\pi]$ is the deviation angle between $v'-v'_*$ and $v-v_*$, and
$B$ is the Boltzmann collision kernel determined by physics 
(related to the cross-section $\Sigma(v-v_*,\sigma)$ 
by the formula $B=|v-v_*| \, \Sigma$). On physical grounds, it is
assumed that $B \geq 0$ and $B$ is a function of $|v-v_*|$ and
$\cos\theta$. 

Boltzmann's collision operator has the fundamental properties of
conserving mass, momentum and energy
  \begin{equation*}
  \int_{\R^N}Q(f,f) \, \phi(v)\,\ud v = 0, \qquad
  \phi(v)=1,v,|v|^2 
  \end{equation*}
and satisfying celebrated Boltzmann's $H$ theorem, which writes formally 
  \begin{equation*} 
  - \dt \int_{\R^N} f \log f \, \ud v = - \int_{\R^N} Q(f,f)\log(f) \, \ud v \geq 0.
  \end{equation*}
The equilibrium distribution is given by the Maxwellian distribution
  \begin{equation*}
  \cal{M}(\rho_\infty,u_\infty,T_\infty)(v)=\frac{\rho_\infty}{(2\pi \, T_\infty)^{N/2}}
  \exp \left( - \frac{\vert u_\infty - v \vert^2} {2 \, T_\infty} \right), 
  \end{equation*}
where $\rho_\infty,\,u_\infty,\,T_\infty$ are the density, mean velocity
and temperature of the gas
  \[
  \rho_\infty = \int_{\tor \times \R^N}f(v) \, \ud x \, \ud v, \quad u_\infty =
  \frac{1}{\rho_\infty}\int_{\tor \times \R^N} v \, f(v) \, \ud x \, \ud v, \] 
  \[ T_\infty = {1\over{N\rho_\infty}}
  \int_{\tor \times \R^N}\vert u_\infty - v \vert^2 \, f(v) \, \ud x \, \ud v, \]
which are determined by the mass, momentum and energy of the initial datum thanks 
to the conservation properties. 

The main physical case of application of this subsection is that of 
hard spheres in dimension $N=3$, where (up to a normalization constant)
 \begin{equation}\label{eq:hs}
 B (|v-v_*|, \cos \theta) = |v-v_*|.
 \end{equation} 
More generally we shall make the following assumption on the collision kernel:
 \begin{itemize}
 \item[{\bf B1.}] We assume that $B$ takes the product form
  \begin{equation}\label{eq:prod}
  B(|v-v_*|,\cos \theta) = \Phi (|v-v_*|) \, b(\cos \theta),
  \end{equation}
 with $\Phi$ and $b$ non-negative and not identically equal to $0$. 
 This decoupling assumption is made for the sake of 
 simplicity and could probably be relaxed at the price of 
 technical complications. 
 \smallskip

 \item[{\bf B2.}]  Concerning the kinetic part, we assume 
 $\Phi$ to be given by 
  \begin{equation}\label{eq:hyprad}
  \Phi (z) = C_\Phi \, z^\gamma
  \end{equation}
 with $\gamma \in [0,1]$. It is customary in physics and mathematics to study the case
 when $\Phi(v-v_*)$ behaves like a power law $|v-v_*|^\gamma$, and one
 traditionally separates between hard potentials ($\gamma >0$),
 Maxwellian potentials ($\gamma = 0$), and soft potentials ($\gamma < 0$). 
 We assume here that we deal with {\bf hard potentials} (or Maxwell molecules). This 
 assumption is crucial since for soft potentials with angular cutoff (see below), 
 the linearized operator  has no spectral gap. 
 \smallskip

 \item[{\bf B3.}] Concerning the angular part, we assume that it is $C^1$ with the controls from above
  \begin{equation}\label{eq:grad}
  \forall \, z \in [-1,1], \ \ \ b(z), \ b'(z) \le C_b.
  \end{equation}
 This implies in particular that $B$ satisfies Grad's {\bf angular cutoff} 
 (see~\cite{Grad58}). 
 Note that the smoothness assumption on $b$ could be relaxed by using truncations and mollifications in the proof.  
 \end{itemize}
 \smallskip

When $b$ is integrable on the sphere $\ens{S}^{N-1}$ (as here 
thanks to {\bf B3}), we define 
  \[ \ell_b := \|b\|_{L^1(\ens{S}^{N-1})} := 
    \big|\ens{S}^{N-2}\big| \, \int_0 ^\pi b(\cos \theta)  \sin^{N-2} \theta \, \ud\theta < +\infty. \]
Without loss of generality we set $\ell_b =1$ in the sequel. 
Then one can split the collision operator in the following way 
  \begin{eqnarray*}
  Q(g,f) &=& Q^+ (g,f) - Q^- (g,f) \\
  Q^+ (g,f) &=& \int _{\R^N \times \ens{S}^{N-1}} \Phi(|v-v_*|) \, b(\cos \theta) \, g'_* f' \, \ud v_* \, \ud\sigma. \\
  Q^-(g,f) &=& \int _{\R^N \times \ens{S}^{N-1}} \Phi(|v-v_*|) \, b(\cos \theta) \, g_* f \, \ud v_* \, \ud\sigma 
  = (\Phi * g) \, f.
  \end{eqnarray*}
We introduce the so-called {\em collision frequency} 
  \begin{equation} \label{eq:colfreq}
  \nu(v) = \int_{\R^N \times \ens{S}^{N-1}} \Phi(|v-v_*|) \, 
  b(\cos \theta) \, \cal{M}(v_*) \, \ud v_* \, \ud\sigma = (\Phi * \cal{M}) (v),
  \end{equation}
and denote by $\nu_0 >0$ the minimum value of $\nu$. 

Using the notation $M = \cal{M}^{1/2}$ the linearized collision operator is given by
    \begin{multline*}
    L(h) = M^{-1} \left[ Q(M h,\cal{M}) + Q(\cal{M},M h) \right] \\
    = M \, \int_{\R^N \times \ens{S}^{N-1}} \Phi(|v-v_*|) \, b(\cos \theta) \, 
    \cal{M}_* \left[ \frac{h' _*}{M'_{*}}+\frac{h'}{M'} - \frac{h_*}{M_{*}} 
        - \frac{h}{M} \right] \, \ud v_* \, \ud\sigma.
    \end{multline*} 
$L$ is self-adjoint on the space $L^2 _v$. 
It splits between a multiplicative part and a non-local part as follows  
 \begin{equation*}
 L(h) = K (h) - \Lambda (h) \ \ \ \mbox{ with } \ \ \ 
 \Lambda (h) = \nu(v) \, h
 \end{equation*}
and 
 \begin{equation*}
 K(h) = L^+ (h) - L^* (h)
 \ \ \ \mbox{ with } \ \ \ 
 L^* (h) = M \left[( h M ) * \Phi \right] 
 \end{equation*}
and
 \begin{equation*}
 L^+ (h) = \int_{\R^N \times \ens{S}^{N-1}} \Phi(|v-v_*|) \, b(\cos \theta) \, 
 \left[ h' M' _{*} + h' _* M' \right] \, M_* \, \ud v_* \, \ud\sigma.
 \end{equation*}
$K$ is bounded and compact in $L^2 _v$, as proved in \cite{Grad63}. 

From the classical spectral theory of $L$ 
it is well-known that with the usual changes of variables
 \begin{multline*}
 \langle L h, h \rangle_{L^2 _{v}}= - \frac{1}{4}\int_{\R^N \times 
 \R^N \times \ens{S}^{N-1}} \Phi(|v-v_*|)\, b(\cos \theta) \, 
  \\ \left[ \frac{h' _*}{M'_{*}}+\frac{h'}{M'} - \frac{h_*}{M_{*}} 
        - \frac{h}{M} \right]^2 \, \cal{M} \, \cal{M}_{*} \, \ud v \, \ud v_{*} \, \ud \sigma \le 0.
 \end{multline*}
This implies that the spectrum of $L$ in $L^2 _{v}$ is included in $\R_-$. 
Moreover the null space of $L$ is 
  \begin{equation}\label{noyauL}
  N(L) = \mbox{Span} \left\{ M, v_1 M,\dots,v_N M, |v|^2 M \right\}.
  \end{equation}

Using the fact that 
 \[ \nu^\Lambda _1 \, (1+|v|)^\gamma \le \nu(v) \le \nu^\Lambda _2 \, (1+|v|)^\gamma \]
for some explicit constants $\nu^\Lambda _1, \nu^\Lambda _2 >0$, and 
that $\nabla_v \nu \in L^\infty$ with explicit bound since $\gamma \in [0,1]$, 
we deduce that assumption {\bf H1} is satisfied with the norm 
 \[ \|h\|_\Lambda = \| h (1+|v|)^{\gamma/2} \|_{L^2}. \] 
Assumption {\bf H1'} can be proved by similar arguments. 

Now we want to prove that 
 \begin{equation} \label{eq:sg}
 \forall \, h \bot N(L), \hspace{0.3cm}
 - \langle h,L h \rangle_{L^2 _{v}} \ge \lambda \, \|h\|^2 _\Lambda. 
 \end{equation}
Controls from below on the collision kernel are necessary to ensure 
the existence of a spectral gap for the linearized operator.  
Under our assumptions, the non-constructive proof of Grad shows that $L$ has 
a spectral gap. Moreover explicit estimates on the spectral gap $\lambda_{L}$ 
have recently been obtained in~\cite{BaMo} and 
extended to explicit estimates of the form~\eqref{eq:sg} in~\cite{Mcoerc}. Following these results
$L$ satisfies {\bf H3} (for the norm $\Lambda$) with explicit bound.  

Now we fix some $\delta>0$ and check that $L$ satisfies assumption {\bf H2} 
(assumption {\bf H2'} can be proved by similar arguments). 

\noindent
Concerning the part $L^*$, this accounts essentially to Young's inequality. 
We easily compute the kernel of the operator
  \begin{equation*}
  L^* h(v) =  \int_{V \in \R^N} h(v+V) \, k^* (v,V) \, \ud V
  \end{equation*}
with
  \begin{equation*}
  k^*(v,V) = M(v)^{1/2} \, \Phi(|V|) \, M(v+V)^{1/2}. 
  \end{equation*}
We introduce the splitting $k^* = k^{*,s} _\var + k^{*,r} _\var$, with 
  \[ k^{*,s} _\var (v,V) = {\cal I}_{\{|V| \ge \var\}} \,  k^* (v,V)  \]
where ${\cal I}$ denotes some mollified indicator function. 
This induces the corresponding decomposition $L^* = L^{*,s} _\var + L^{*,r} _\var$. 
It is straightforward that 
  \[ \| L^{*,r} _\var \|_{L^2 \to L^2} \xrightarrow[]{\var \to 0} 0 \]
and
  \[ \left| \nabla_v k^{*,s} _\var \right|, \ \ 
     \left| \nabla_V k^{*,s} _\var \right|, \ \ 
     \left| \nabla_v k^{*,r} _\var \right| \le C(\var) \,  M(V)^{1/8}. \]
Hence we deduce 
  \begin{equation}\label{eq:regL^*1}
  \left\| \nabla_{v} L^{*,s} _\var h \right\|_{L^2} \le C(\var) \, \|h\|_{L^2} 
  \end{equation}
and 
  \begin{equation}\label{eq:regL^*2}
  \left\| \nabla_{v} L^{*,r} _\var h \right\|_{L^2} \le \delta \, \|h\|_{H^1} + 
                                    C(\var) \, \|h\|_{L^2} 
  \end{equation}
if $\var$ is small enough.

Now we turn to the part $L^+$.  
We follow Grad computations~\cite[Sections 2 and 3]{Grad63} 
(recalled also in~\cite[Chapter 7, Section 2]{CIP94}) to compute the kernel 
of $L^+$, and apply the same kind of estimates as in \cite{Meepts}. 
We make the changes the variables 
\begin{itemize}
\item $\sigma \in \ens{S}^{N-1}$, $v_* \in \R^N$ $\longrightarrow$ $\omega = (v'-v)/|v'-v| \in \ens{S}^{N-1}$, 
$v_* \in \R^N$: the jacobian amounts to change $b$ into 
  \[ \tilde b (\theta) = 2^{N-1} \, \sin^{N-2} \theta/2 \, b(\theta); \] 
\item then $\omega \in \ens{S}^{N-1}$, $v_* \in \R^N$ $\longrightarrow$ 
$\omega \in \ens{S}^{N-1}$, $u=v-v_* \in \R^N$:   
the jacobian is equal to $1$; 
\item then keeping $\omega$ fixed, decompose orthogonally $u = u_0 \omega +W$ with 
$u_0 \in \R$ and $W \in \omega^\bot$: the jacobian is equal to $1$; 
\item finally keeping $W \in V^\bot$ fixed, $\omega \in \ens{S}^{N-1}$, $u_0 \in \R$ 
$\longrightarrow$ $V = u_0 \omega \in \R^N$: the jacobian is $(1/2)|V|^{-(N-1)}$.
\end{itemize} 
Thus we get
  \begin{equation*}
  L^+ h(v) =  \int_{V \in \R^N} h(v+V) \, k^+ (v,V) \, \ud V
  \end{equation*}
with
  \begin{multline*}
  k^+(v,V) = \mbox{cst} \, |V|^{-(N-1)} \,  \int_{W \in V^\bot} \Phi(\sqrt{|V|^2 + |W|^2}) \, 
        \tilde b\left( \frac{|W|^2-|V|^2}{|W|^2+|V|^2} \right) \times \\ 
        M(v+W)^{1/2} \,  M(v+V+W)^{1/2} \, \ud W. 
  \end{multline*}
This kernel can be written as 
  \begin{multline*}
  k^+(v,V) = \mbox{cst} \, M(V)^{1/4} \, |V|^{-(N-1)} \, 
        \int_{W \in V^\bot} \Phi(\sqrt{|V|^2 + |W|^2}) \, 
        \tilde b\left( \frac{|W|^2-|V|^2}{|W|^2+|V|^2} \right) \times \\ 
        M(v+V+W/2) \, \ud W. 
  \end{multline*}
Moreover it is shown in~\cite[Chapter 7, Section 2]{CIP94} that 
  \[ \| {\bf 1}_{|\cdot|\ge R} L^+ \|_{L^2 \to L^2} \xrightarrow[]{R \to \infty} 0. \]
We use this to perform the splitting 
  \[ k^+ = k^{+,s} _\var + k^{+,r} _\var \]
with 
  \[ k^{+,s} _\var (v,V) = {\cal I}_{\{|v| \le \var^{-1}\}} \, {\cal I}_{\{|V| \ge \var \}} 
                            \,  k^+ (v,V),  \]
where ${\cal I}$ denotes some mollified indicator function. The corresponding 
decomposition of $L$ is denoted by
  \[ L^+ = L^{+,s} _\var + L^{+,r} _\var. \]
It is straightforward that 
  \[ \| L^{+,r} _\var \|_{L^2 \to L^2} \xrightarrow[]{\var \to 0} 0 \]
and 
  \[ \left| \nabla_v k^{+,s} _\var \right|, \ \ 
     \left| \nabla_V k^{+,s} _\var \right|, \ \ 
     \left| \nabla_v k^{+,r} _\var \right| \le C(\var) \,  M(V)^{1/8}. \]
Hence we deduce 
  \begin{equation}\label{eq:regL^+1}
  \left\| \nabla_{v} L^{+,s} _\var h \right\|_{L^2} \le C(\var) \, \|h\|_{L^2} 
  \end{equation}
and 
  \begin{equation}\label{eq:regL^+2}
  \left\| \nabla_{v} L^{+,r} _\var h \right\|_{L^2} \le \delta \, \|h\|_{H^1} + 
                                    C(\var) \, \|h\|_{L^2} 
  \end{equation}
as long as $\var$ is small enough.\\
This concludes the proof by gathering~\eqref{eq:regL^*1}, \eqref{eq:regL^*2}, \eqref{eq:regL^+1} 
and \eqref{eq:regL^+2}. 

Finally let us consider the bilinear part given by
  \begin{multline*}
  \Gamma(h_1,h_2) = M^{-1} \left[ Q(M h,M h) + Q(M h,M h) \right] \\
    = \int_{\R^N \times \ens{S}^{N-1}} \Phi(|v-v_*|) \, b(\cos \theta) \, 
       M_* \left[ (h_1)' _* (h_2)' - (h_1)_* (h_2)  \right] \, \ud v_* \, \ud\sigma \\
    = \int_{\R^N \times \ens{S}^{N-1}} \Phi(|u|) \, b(\cos \theta) \, 
       M_* \left[ (h_1)' _* (h_2)' - (h_1)_* (h_2)  \right] \, \ud u \, \ud\sigma
    \end{multline*} 
with the notation $u=v-v_*$. We estimate 
  \begin{multline*}
    \int_{\R^N} \Gamma(h_1,h_2) \, \varphi \, \ud v \\
    \le C \, \|\varphi\|_{L^2} \, 
        \left( \int_{\R^N} \left| \int_{\R^N \times \ens{S}^{N-1}} |(h_1)'| \, |(h_2)'_*| \, M(v+u) \, 
               \langle u \rangle^\gamma \, \ud u \, \ud \sigma \right|^2 \, \ud v \right)^{1/2} \\
    \le C \, \|\varphi\|_{L^2} \, 
        \left( \int_{\R^N \times \R^N \times \ens{S}^{N-1}} |(h_1)'|^2 \, |(h_2)'_*|^2 \, 
               \langle v \rangle^\gamma \, \ud v \, \ud u \, \ud \sigma \right) \\
    \le C \, \|\varphi\|_{L^2} \, 
        \left( \|h_1\|_{L^2 _v} \, \|h_2\|_{\Lambda_v} + \|h_1\|_{\Lambda_v} \, \|h_2\|_{L^2_v} \right), 
  \end{multline*} 
which implies 
  \[ \| \Gamma(h_1,h_2) \|_{L^2 _v} \le C \, 
     \left( \|h_1\|_{L^2 _v} \, \|h_2\|_{\Lambda_v} + \|h_1\|_{\Lambda_v} \, \|h_2\|_{L^2_v} \right). \]
Together with Leibnitz formula to differentiate $\Gamma$ according to $v$ and $x$ and Sobolev 
embeddings this concludes the proof of {\bf H4} for $E(k_0/2) > N/2$. 

\subsection{The Landau equation}\label{ss:La}

This subsection deals with the Landau equation (for $N \ge 2$)
  \begin{equation}\label{eq:Lan}
  \partial_{t} f + v \cdot \nabla_{x} f = Q(f,f), \qquad t \ge 0,\ x \in \tor, \ v \in \R^N
  \end{equation} 
which features the collision operator
  \[ Q (f,f)(v) = \nabla _v \cdot \left( \int_{\R^N} 
  {\bf A}(v-v_*) \left[ f_* \left( \nabla f \right) 
  - f \left( \nabla f \right)_* \right] \, \ud v_* \right),
  \]
where ${\bf A} (z) = |z|^2 \, \Phi(|z|) \, {\bf P}(z)$, $\Phi$ is a non-negative function, 
and ${\bf P}(z)$ is the orthogonal projection onto $z^\bot$, {\it i.e.},  
 \begin{equation*}
 \left( {\bf P}(z) \right) _{i,j} = \delta_{i,j} - \frac{z_i z_j}{|z|^2}
 \end{equation*}
We use again the notation $f_{*}= f(v_{*})$. 
This operator is used for instance in models for plasma. 
In this case the interaction among the particles is {\it via} the Coulomb 
potential and $\Phi(|z|) = |z|^{-3}$ in dimension $3$. 
For more details see~\cite[Chapter~1, Section~1.7]{Vi:hand} and the references therein. 
Indeed in this case the Boltzmann collision operator 
does not make sense anymore (see~\cite[Annex~I, Appendix]{Vi:habil}). 

Landau's collision operator has the fundamental properties of
conserving mass, momentum and energy
  \begin{equation*}
  \int_{\R^N}Q(f,f) \, \phi(v)\,\ud v = 0, \qquad
  \phi(v)=1,v,|v|^2 
  \end{equation*}
and satisfying Boltzmann's $H$ theorem, which writes formally 
  \begin{equation*} 
  - \dt \int_{\R^N} f \log f \, \ud v = - \int_{\R^N} Q(f,f)\log(f) \, \ud v \geq 0.
  \end{equation*}
The equilibrium distribution is given by the Maxwellian distribution
  \begin{equation*}
  \cal{M}(\rho_\infty,u_\infty,T_\infty)(v)=\frac{\rho_\infty}{(2\pi \, T_\infty)^{N/2}}
  \exp \left( - \frac{\vert u_\infty - v \vert^2} {2 \, T_\infty} \right), 
  \end{equation*}
where $\rho_\infty,\,u_\infty,\,T_\infty$ are determined as in the Boltzmann case. 
\smallskip

We make the following assumption on the collision kernel:
 \begin{itemize}
 \item[{\bf L1.}]  We assume $\Phi$ to be given by 
  \begin{equation}\label{eq:Landhyprad}
  \Phi (z) = C_\Phi \, z^\gamma
  \end{equation}
 with $\gamma \in [-2,1]$. By analogy with the Boltzmann equation, one could 
 say that this assumption covers hard and moderately soft potentials. 
 \end{itemize}
 \smallskip

We consider fluctuations around equilibrium of the form 
$f=\cal{M} + M h$. The linearized collision operator is given by
    \begin{equation*}
    L (h)= M^{-1} \, \nabla_v \cdot \left( \int_{\R^N} 
    {\bf A}(v-v_*) \left[ \left( \frac{\nabla_v h}M \right)  - \left( \frac{\nabla_v h}M \right)_* \right] 
    \, \cal{M} \cal{M}_* \, \ud v_* \right).
    \end{equation*} 
$L$ is self-adjoint on the space $L^2 _v$. 
It splits between an (almost) ``convolution part'' and a diffusive part: 
 \begin{equation*}
 L(h) = K (h) - \Lambda (h)
 \end{equation*}
with 
 \begin{equation*}
 K(h) = - M^{-1} \, \nabla_v \cdot \left( \int_{\R^N} 
    {\bf A}(v-v_*) \left( \frac{\nabla_v h}M \right)_* \, \cal{M} \cal{M}_* \, \ud v_* \right),
 \end{equation*}
and
 \begin{equation*}
 \Lambda(h) = - M^{-1} \, \nabla_v \cdot \left( \int_{\R^N} 
    {\bf A}(v-v_*) \left( \frac{\nabla_v h}M \right) \, \cal{M} \cal{M}_* \, \ud v_* \right).
 \end{equation*}
Estimate {\bf H2} on $K$ is easily verified since 
 \begin{equation*}
 K(h) = \int_{\R^N} k(v,v_*) \, h_* \, \ud v_* 
 \end{equation*}
where the kernel 
 \[ k(v,v_*) = \left[ \left(\nabla_v\right)^T \left(\frac{{\bf A}(v-v_*) \cal{M} \cal{M}_*}{M M_*}\right) 
  \left(\nabla_v \right) \right] \]
belongs straightforwardly to $L^2(\R^N \times \R^N)$ and also to 
$H^1 (\R^N \times \R^N)$ except possibly for a small 
region $(v-v_*) \sim 0$ which can split as in the Boltzmann case. 

It is well-known from the classical spectral theory 
of $L$ that with the usual changes of variables we have
 \begin{multline*}
 \langle h,L h \rangle_{L^2 _{v}}= - \frac{1}{2}\int_{\R^N \times \R^N} 
 \Phi(|v-v_*|)\, |v-v_*|^2 \, 
  \\ \left\| {\bf P} \left[ \left( \frac{\nabla_v h}M \right) 
  - \left( \frac{\nabla_v h}M \right)_* \right] \right\| ^2 
  \, \cal{M} \, \cal{M}_{*} \, \ud v \, \ud v_{*} \le 0.
 \end{multline*}
This implies that the spectrum of $L$ in $L^2 _{v}$ is included in $\R_-$. 
Moreover the null space of $L$ is 
  \begin{equation}\label{noyauLFP}
  N(L) = \mbox{Span} \left\{ M, v_1 M,\dots,v_N M, |v|^2 M \right\}.
  \end{equation}

Now let us use some estimates proved in~\cite{Guo:LE:02}. 
First we define the norm
  \[ \| h \|_{\Lambda_v} ^2 = \Big\| h \langle v \rangle^{1+\gamma/2} \Big\|_{L^2 _v} ^2 
         + \Big\| \big({\bf P}(v) \nabla_v h\big) \langle v \rangle^{1+\gamma/2} \Big\|_{L^2 _v} ^2 
         +  \Big\| \big( (1- {\bf P}(v)) \nabla_v h\big) \langle v \rangle^{\gamma/2} \Big\|_{L^2 _v} ^2 \]
which is stronger than $L^2 _v$ as soon as $\gamma \ge -2$. 
In~\cite[Section~2]{Guo:LE:02} it is proven that
  \[ \langle \Lambda h, h \rangle_{L^2 _v} \ge C \,  \| h \|_{\Lambda_v} ^2 \]
with explicit constant, from which {\bf H1} follows. 
For the bilinear term Theorem~3 from~\cite{Guo:LE:02} 
together with Sobolev embeddings yields {\bf H4} in the norm $\Lambda$ with explicit 
constant as long as $E(k_0/2) > N/2$. Again the stronger assumptions {\bf H1'-H2'} 
are deduced straightforwardly with the same arguments. 
Finally assumption {\bf H3} is proved in~\cite[Section~2, Lemma~5]{Guo:LE:02} by 
non-constructive arguments (and an explicit proof is given in~\cite{MS**}). 
This concludes the proof. 


\subsection{Remarks on other models} 
Linear models of radiative transfer in the torus enter straightforwardly our 
abstract framework. It is likely that  
linear scattering Boltzmann models or semi-conductors collisional models also do so. 
Moreover it is easy to see on the linear relaxation models (as well as on more 
general linear scattering models) that one could add with very few changes in our proof 
some scattering rate $\Sigma=\Sigma(x)$ depending on $x$ in front of the collision operator:  
assuming that $\Sigma \in C^\infty$ and 
  \[ \forall \, x \in \tor, \quad 0< \Sigma_- \le \Sigma(x) \le \Sigma_+ < +\infty \]
for some constants $\Sigma_-, \Sigma_+ >0$, the conclusion of Theorem~\ref{theo:T} still holds.  
\medskip

For the Boltzmann equation with soft potentials and Grad's angular cutoff, smooth 
solutions near the Maxwellian equilibrium have been built in~\cite{Guo:BE:03}: 
by including polynomial weight in $v$ depending on the order of 
the derivatives in the energy estimates, it is likely that 
one can adapt our proof to build a norm which is decreasing along the flow inspiring from~\cite{Guo:BE:03}. 
However in this case the integro-differential operator 
$T$ is not coercive for this norm, instead it satisfies degenerated coercivity estimates 
for some weaker norms. This is enough to built smooth solution, but does not yield exponential 
convergence towards equilibrium. Nevertheless as noticed in~\cite{SG} one can deduce from it 
polynomial rates of decay to equilibrium by interpolating between a ladder of norms. 
\medskip

Our analysis works at the linear level for the linearized Boltzmann equation for hard potentials 
without Grad's angular cutoff assumption, using explicit spectral gap estimates on $L$ provided 
by~\cite{BaMo} (note also that it could cover some moderately soft potentials interactions without 
Grad's angular cutoff, using explicit spectral gap estimates on $L$ provided by~\cite{MS**}). 
Indeed in~\cite{Mcoerc} (see also~\cite{MS**}), it is shown how to write the linearized collision 
operator in the form $K-\Lambda$ with some regularizing $K$ and some 
coercive $\Lambda$, and how to obtain coercivity estimates on $L$. 
However at now it is not known how to control the non-linear term in terms 
of a coercivity norm $\Lambda$ adapted to the linearized 
operator. The functional space of these coercivity estimates in~\cite{Mcoerc} 
is a {\em local} Sobolev space with the right fractional order, 
but which does not seem sufficient to control the non-linear term.   
\medskip

\bigskip
\noindent
{\bf{Acknowledgment}}: Support by the European network HYKE, funded by the EC as
contract HPRN-CT-2002-00282, is acknowledged. The second author acknowledges 
support by the FWF-Doktoratskolleg ``Differential Equation Models in Science 
and Engineering'' and travelling support by the \"{O}AD (Austrian French Cooperation, 
Amad\'{e}e Project No. 19/2003). Authors would like to thank Mar\'{\i}a Jos\'e 
C\'aceres and C\'edric Villani for fruitful discussions. In particular they thank 
C\'edric Villani for pointing out and discussing his notes~\cite{Vi:hypofd,Vi:hypo}, from which 
the idea to use time derivatives of some mixed terms in the energy estimates originates.     

\begin{flushleft} \signcm \end{flushleft} 
\begin{flushleft} \signln \end{flushleft} 
\end{document}